\documentclass{amsart}
\usepackage{graphicx}
\usepackage{amscd}
\usepackage{amsmath}
\usepackage{amsfonts}
\usepackage{amssymb}
\newtheorem{theorem}{Theorem}
\theoremstyle{plain}

\newtheorem{corollary}{Corollary}

\newtheorem{definition}{Definition}

\newtheorem{lemma}{Lemma}

\newtheorem{proposition}{Proposition}

\numberwithin{equation}{section}

\begin{document}
\title{Symmetric Functions and $B_{N}$-Invariant Spherical Harmonics}
\author{Charles F. Dunkl}
\address{Department of Mathematics, University of Virginia, PO Box 400137,
Charlottesville VA 22904-4137}
\email{cfd5z@virginia.edu}
\urladdr{http://www.math.virginia.edu/\symbol{126}cfd5z/}
\thanks{During the preparation of this paper the author was partially supported by NSF
grants DMS 9970389 and DMS 0100539.}
\subjclass{33C55, 05E05, 81Q05, 82D25}
\date{July 15, 2002}

\begin{abstract}
The wave functions of a quantum isotropic harmonic oscillator in N-space
modified by barriers at the coordinate hyperplanes can be expressed in terms
of certain generalized spherical harmonics. These are associated with a
product-type weight function on the sphere. Their analysis is carried out by
means of differential-difference operators. The symmetries of this system
involve the Weyl group of type B, generated by permutations and changes of
sign of the coordinates. A new basis for symmetric functions as well as an
explicit transition matrix to the monomial basis is constructed. This basis
leads to a basis for invariant spherical harmonics. The determinant of the
Gram matrix for the basis in the natural inner product over the sphere is
evaluated, and there is a formula for the evaluation of the basis elements at
(1,1,...,1). When the underlying parameter is specialized to zero, the basis
consists of ordinary spherical harmonics with cube group symmetry, as used for
wave functions of electrons in crystals. The harmonic oscillator can also be
considered as a degenerate interaction-free spin Calogero model.
\end{abstract}\maketitle

\section{Introduction}

There are interesting families of potentials invariant under permutations and
changing of signs of coordinates. The most basic one is a central potential
perturbed by a crystal field with cubic symmetry. Another important example is
the spin Calogero-Moser system of $B$-type. In this paper we study a potential
which can be considered as an $N$-dimensional isotropic harmonic oscillator
modified by barriers at the coordinate hyperplanes, or as a degenerate
Calogero-Moser model with no interaction. The main object is to study
invariant harmonic polynomials, which when multiplied by radial Laguerre
polynomials provide a complete decomposition of the invariant wave functions.
This is made possible by use of the author's differential-difference operators
and the construction of a new basis for symmetric functions.

We consider generalized spherical harmonic polynomials invariant under the
action of the hyperoctahedral group $B_{N}$ acting on $\mathbb{R}^{N}$. Let
$\mathbb{N}_{0}=\left\{  0,1,2,\ldots\right\}  $, the set of compositions with
$N$ parts is $\mathbb{N}_{0}^{N}$; for $\alpha\in\mathbb{N}_{0}^{N}$ let
$\left|  \alpha\right|  =\sum_{i=1}^{N}\alpha_{i}$, $\alpha!=\prod_{i=1}%
^{N}\alpha_{i}!$ and $\left(  t\right)  _{\alpha}=\prod_{i=1}^{N}\left(
t\right)  _{\alpha_{i}}$, (the Pochhammer symbol is $\left(  t\right)
_{n}=\prod_{i=1}^{n}\left(  t+i-1\right)  $). Denote the cardinality of a
finite set $E$ by $\#E$. The set $\mathcal{P}$ of partitions consists of
finite sequences $\lambda=\left(  \lambda_{1},\lambda_{2},\ldots\right)
\in\bigcup_{N=1}^{\infty}\mathbb{N}_{0}^{N}$ satisfying $\lambda_{1}%
\geq\lambda_{2}\geq\ldots\geq0$ and two partitions with the same nonzero
components are identified. Then $\mathcal{P}_{n}=\left\{  \lambda
\in\mathcal{P}:\left|  \lambda\right|  =n\right\}  $ for $n=1,2,3,\ldots$. For
a partition the length $l\left(  \lambda\right)  =\#\left\{  i:\lambda_{i}%
\geq1\right\}  $ is the number of nonzero components. For a given $N$ we will
need the partitions with $l\left(  \lambda\right)  \leq N$, denoted by
$\mathcal{P}^{\left(  N\right)  }=\mathcal{P\cap}\mathbb{N}_{0}^{N},$ and
$\mathcal{P}_{n}^{\left(  N\right)  }=\mathcal{P}_{n}\mathcal{\cap}%
\mathbb{N}_{0}^{N}$. One partial ordering for partitions is given by
containment of Ferrers diagrams: thus $\mu\subset\lambda$ means $\mu_{i}%
\leq\lambda_{i}$ for each $i$. The notation $\varepsilon_{i}=\left(
0,\ldots,\overset{i}{1},\ldots\right)  $ for the $i^{th}$ standard basis
vector in $\mathbb{N}_{0}^{N}$ (for $1\leq i\leq N$) is convenient for
describing contiguous partitions; for example $\lambda+\varepsilon_{1}=\left(
\lambda_{1}+1,\lambda_{2},\ldots\right)  $.

For $x\in\mathbb{R}^{N}$ and $\alpha\in\mathbb{N}_{0}^{N}$ let $x^{\alpha
}=\prod_{i=1}^{N}x_{i}^{\alpha_{i}}$, a monomial of degree $\left|
\alpha\right|  $. Let $\mathbb{P}_{n}^{\left(  N\right)  }=\mathrm{span}%
\left\{  x^{\alpha}:\alpha\in\mathbb{N}_{0}^{N},\left|  \alpha\right|
=n\right\}  $, the space of homogeneous polynomials of degree $n\geq0$ in $N$
variables. The Laplacian is $\Delta=\sum_{i=1}^{N}\frac{\partial^{2}}{\partial
x_{i}^{2}}$, and the Euclidean norm is $\left|  \left|  x\right|  \right|
=\left(  \sum_{i=1}^{N}x_{i}^{2}\right)  ^{1/2}$. Fix the parameter
$\kappa\geq0.$

We will be concerned with the operator $\Delta+2\kappa\sum_{i=1}^{N}\frac
{1}{x_{i}}\frac{\partial}{\partial x_{i}}$. It is associated with the
Calogero-Sutherland model with the potential function $V\left(  x\right)
=\omega^{2}\left|  \left|  x\right|  \right|  ^{2}+\kappa\left(
\kappa-1\right)  \sum\limits_{i=1}^{N}\dfrac{1}{x_{i}^{2}}$ (with $\omega>0$).
This is an $N$-dimensional isotropic harmonic oscillator modified by barriers
at the coordinate hyperplanes; or it could be considered as a degenerate form
of the type-$B$ spin model for $N$ particles with no interaction. The spin
model with interactions and reflecting barriers was studied by Yamamoto and
Tsuchiya \cite{YT}. With the base state
\[
\psi\left(  x\right)  =\exp\left(  -\frac{\omega}{2}\left|  \left|  x\right|
\right|  ^{2}\right)  \prod_{i=1}^{N}\left|  x_{i}\right|  ^{\kappa}%
\]
we obtain the conjugate of the Hamiltonian, for any smooth function $f$ on
$\mathbb{R}^{N}$
\begin{align*}
&  \psi\left(  x\right)  ^{-1}\left(  -\Delta+\omega^{2}\left|  \left|
x\right|  \right|  ^{2}+\kappa\left(  \kappa-1\right)  \sum_{i=1}^{N}\frac
{1}{x_{i}^{2}}\right)  \left(  \psi\left(  x\right)  f\left(  x\right)
\right) \\
&  =\left(  -\Delta-2\kappa\sum_{i=1}^{N}\frac{1}{x_{i}}\frac{\partial
}{\partial x_{i}}+N\omega\left(  2\kappa+1\right)  +2\omega\sum_{i=1}^{N}%
x_{i}\frac{\partial}{\partial x_{i}}\right)  f\left(  x\right)  .
\end{align*}

The $\mathbb{Z}_{2}^{N}$ version of the differential-difference operators
introduced by the author (\cite{D1}, or see \cite[Ch. 4]{DX}) are the
fundamental tools for analyzing this eigenfunction problem. The action of the
group $B_{N}$ on $\mathbb{R}^{N}$ induces an action on functions, denoted by
$wf\left(  x\right)  =f\left(  xw\right)  $ for $w\in B_{N}.$ For $1\leq i\leq
N$ let $\sigma_{i}$ denote the reflection
\[
\left(  x_{1},\ldots,x_{N}\right)  \sigma_{i}=\left(  x_{1},\ldots
,-x_{i},\cdots,x_{N}\right)  ,
\]
and define the first-order operator $\mathcal{D}_{i}$ by
\[
\mathcal{D}_{i}f\left(  x\right)  =\frac{\partial}{\partial x_{i}}f\left(
x\right)  +\kappa\frac{f\left(  x\right)  -f\left(  x\sigma_{i}\right)
}{x_{i}}%
\]
for sufficiently smooth functions $f$ on $\mathbb{R}^{N}$. The Laplacian
operator is
\[
\Delta_{\kappa}=\sum_{i=1}^{N}\mathcal{D}_{i}^{2},
\]
then
\[
\Delta_{\kappa}f\left(  x\right)  =\Delta f\left(  x\right)  +\kappa\sum
_{i=1}^{N}\left(  \frac{2}{x_{i}}\frac{\partial}{\partial x_{i}}f\left(
x\right)  -\frac{f\left(  x\right)  -f\left(  x\sigma_{i}\right)  }{x_{i}^{2}%
}\right)  .
\]
If $f$ is even in each $x_{i}$ (that is, $\mathbb{Z}_{2}^{N}$-invariant) then
$\Delta_{\kappa}f\left(  x\right)  =\Delta f\left(  x\right)  +$%
\newline $2\kappa\sum_{i=1}^{N}\frac{1}{x_{i}}\frac{\partial}{\partial x_{i}%
}f\left(  x\right)  $. If we modify the above Hamiltonian to
\[
\mathcal{H}=-\Delta+\omega^{2}\left|  \left|  x\right|  \right|  ^{2}%
+\kappa\sum_{i=1}^{N}\frac{\kappa-\sigma_{i}}{x_{i}^{2}},
\]
which is identical to the previous one when acting on $\mathbb{Z}_{2}^{N}%
$-invariant functions, we obtain
\[
\psi^{-1}\mathcal{H}\psi f\left(  x\right)  =\left(  -\Delta_{\kappa}%
+2\omega\sum_{i=1}^{N}x_{i}\frac{\partial}{\partial x_{i}}+\omega N\left(
2\kappa+1\right)  \right)  f\left(  x\right)  .
\]
The eigenfunctions of this operator can be expressed as products of radial
Laguerre polynomials and homogeneous harmonic polynomials.

\begin{definition}
For $n\geq0$ the space of harmonic (sometimes called $h$-harmonic) homogeneous
polynomials of degree $n$ is $\mathbb{H}_{n}=\left\{  f\in\mathbb{P}%
_{n}^{\left(  N\right)  }:\Delta_{\kappa}f=0\right\}  .$ Let $\mathbb{H}%
_{2n}^{0}=\left\{  f\in\mathbb{H}_{2n}:\sigma_{i}f=f,1\leq i\leq N\right\}  $;
this is the space of $\mathbb{Z}_{2}^{N}$-invariant harmonic polynomials, even
in each $x_{i},1\leq i\leq N$. The subspace of $\mathbb{H}_{2n}^{0}$
consisting of polynomials invariant under permutation of coordinates is
denoted by $\mathbb{H}_{2n}^{B}$.
\end{definition}

It is not hard to show that if $f\in\mathbb{H}_{n}$ then $L_{s}^{\left(
n+N/2+N\kappa-1\right)  }\left(  \omega\left|  \left|  x\right|  \right|
^{2}\right)  f\left(  x\right)  $ is an eigenfunction of $\psi^{-1}%
\mathcal{H}\psi$ with eigenvalue $\omega\left(  2n+4s+N\left(  2\kappa
+1\right)  \right)  $ (this calculation uses equation \ref{dpower}). So
combining bases for each $\mathbb{H}_{n}$ with Laguerre polynomials provides a
basis for the polynomial eigenfunctions (positive energy) of $\psi
^{-1}\mathcal{H}\psi$. The use of functions of this type in the general
Calogero-Moser models of types $A$ and $B$ is discussed in van Diejen
\cite{vD}.

The classical motion problem corresponding to $\mathcal{H}$ is easily solved:
in the one-dimensional case the particle with mass 1 at $s\in\mathbb{R}$
satisfies $\left(  \frac{d}{dt}\right)  ^{2}s\left(  t\right)  =-\frac
{\partial}{\partial s}V\left(  s\right)  =-2\omega^{2}s+2\kappa\left(
\kappa-1\right)  s^{-3}$ with the solution (for $s>0,\kappa>1$)
\begin{align*}
s\left(  t\right)   &  =\left(  q+a\sin\left(  2^{3/2}\omega\left(
t-t_{0}\right)  \right)  \right)  ^{1/2}\\
q  &  =\left(  a^{2}+\kappa\left(  \kappa-1\right)  /\omega^{2}\right)
^{1/2},
\end{align*}
here $t_{0}$ is an arbitrary phase shift, $a\geq0$ is arbitrary, and the
energy is $\frac{1}{2}\left(  \frac{ds}{dt}\right)  ^{2}+V\left(  s\right)
=2q\omega^{2}$.

Primarily we will study the invariant functions associated with $\Delta
_{\kappa}$; the invariance is with respect to the Weyl group of type
\textit{B, }thus the functions are invariant under sign changes $\left\{
\sigma_{i}\right\}  $ and permutations of coordinates. Such functions are
expressed as symmetric functions of the variable $x^{2}=\left(  x_{1}%
^{2},\ldots,x_{N}^{2}\right)  $. Note that a special case of this study is the
problem of spherical harmonics on $\mathbb{R}^{3}$ which satisfy $B_{3}%
$-invariance, appearing in the wave functions of electrons in crystals. We
will construct an explicit basis for these polynomials, but unfortunately it
is not orthogonal. There seems to be a good reason why an orthogonal basis has
not yet been found: let $R_{ij}=x_{i}\mathcal{D}_{j}-x_{j}\mathcal{D}_{i}$
(then $R_{ij}\Delta_{\kappa}=\Delta_{\kappa}R_{ij}$, and $\sqrt{-1}\,R_{ij}$
is an angular momentum operator), and let $\mathcal{S}_{n}=\sum_{1\leq i<j\leq
N}R_{ij}^{2n}$, then $\mathcal{S}_{1}$ is the Casimir (or Laplace-Beltrami)
operator, for each $n\ $the operator $\mathcal{S}_{n}$ is self-adjoint in the
natural inner product on the sphere (defined later) but $\mathcal{S}%
_{2},\mathcal{S}_{3}$ do not commute, already for $N=3 $. So the usual
machinery for constructing good orthogonal bases (like Jack polynomials) does
not work here. The conjugate
\[
\psi R_{ij}\psi^{-1}=x_{i}\frac{\partial}{\partial x_{j}}-x_{j}\frac{\partial
}{\partial x_{i}}-\kappa\left(  \frac{x_{i}}{x_{j}}\sigma_{j}-\frac{x_{j}%
}{x_{i}}\sigma_{i}\right)
\]
commutes with the Hamiltonian $\mathcal{H}$ (see Taniguchi \cite{kT} for a
more general treatment of commuting operators in the context of $r^{-2}$ type
potentials). By straightforward calculation
\begin{align*}
\mathcal{S}_{1}  &  =\sum_{1\leq i<j\leq N}R_{ij}^{2}\\
&  =\left|  \left|  x\right|  \right|  ^{2}\Delta_{\kappa}-\left(  \sum
_{i=1}^{N}x_{i}\mathcal{D}_{i}\right)  ^{2}-\left(  N-2+2\kappa\sum_{i=1}%
^{N}\sigma_{i}\right)  \left(  \sum_{i=1}^{N}x_{i}\mathcal{D}_{i}\right)  .
\end{align*}
Thus $\mathcal{S}_{1}$ has the eigenvalue $-2n\left(  N-2+2n+2N\kappa\right)
$ on $\mathbb{H}_{2n}^{0}$.

\section{A basis for symmetric functions}

In this section the number $N$ of variables is not specified, with the
understanding that it is not less than the length of any partition that
appears. We use the notation of Macdonald \cite{M} for symmetric polynomials
in the variable $x\in\mathbb{R}^{N}$. Let $S_{N}$ denote the symmetric group
on $N$ objects, considered as the group of $N\times N$ permutation matrices,
acting on the left on $\mathbb{N}_{0}^{N}$. For any $\alpha\in\mathbb{N}%
_{0}^{N}$ let $\alpha^{+}$ denote the unique partition $w\alpha\in
\mathcal{P}^{\left(  N\right)  }$, for some $w\in S_{N}$ (the sorting of the
components of $\alpha$ in nonincreasing order; $w$ need not be unique).

\begin{definition}
For $\lambda\in\mathcal{P}^{\left(  N\right)  }$ the monomial symmetric
function is
\[
m_{\lambda}=\sum\left\{  x^{\alpha}:\alpha\in\mathbb{N}_{0}^{N},\alpha
^{+}=\lambda\right\}  ,
\]
summing over all distinct permutations of $\lambda$. The elementary symmetric
function of degree $1$ is
\[
e_{1}=\sum_{i=1}^{N}x_{i}.
\]
\end{definition}

It turns out that the basis elements for invariant harmonics are labeled by
partitions with the property $\lambda_{1}=\lambda_{2}$; that is,
$\dim\mathbb{H}_{2n}^{B}=\allowbreak\#\left\{  \lambda\in\mathcal{P}%
_{n}^{\left(  N\right)  }:\lambda_{1}=\lambda_{2}\right\}  $. Further the
formula for the projection onto harmonics uses powers of $\left|  \left|
x\right|  \right|  ^{2}=e_{1}\left(  x_{1}^{2},\ldots,x_{N}^{2}\right)  .$
This leads to the following definition of a basis for symmetric functions well
suited for the present study.

\begin{definition}
For $\lambda\in\mathcal{P}^{\left(  N\right)  }$ the modified monomial
symmetric function is
\[
\widetilde{m}_{\lambda}=e_{1}^{\lambda_{1}-\lambda_{2}}m_{\left(  \lambda
_{2},\lambda_{2},\lambda_{3},\ldots\right)  }.
\]
\end{definition}

We will show that $\left\{  \widetilde{m}_{\lambda}:\lambda\in\mathcal{P}%
^{\left(  N\right)  }\right\}  $ is a basis for the symmetric polynomials on
$\mathbb{R}^{N}$ and the transition matrix to the $m$-basis has entries in
$\mathbb{Z}$, is unimodular and triangular in a certain ordering. One
direction is an easy consequence of the multinomial theorem. We use the
notation $\binom{j}{\alpha}$ for the multinomial coefficient, where $\alpha
\in\mathbb{Z}^{N},\sum_{i=1}^{N}\alpha_{i}=j$, and $\binom{j}{\alpha}%
=\frac{j!}{\alpha!}$ if each $\alpha_{i}\geq0$, else $\binom{j}{\alpha}=0$.

\begin{definition}
\label{e1m}For $\lambda,\nu\in\mathcal{P}^{\left(  N\right)  }$ the
coefficient $%
\genfrac{\langle}{\rangle}{0pt}{}{\lambda}{\nu}%
$ is defined by the expansion
\[
e_{1}^{j}m_{\nu}=\sum_{\left|  \lambda\right|  =\left|  \nu\right|  +j}%
\genfrac{\langle}{\rangle}{0pt}{}{\lambda}{\nu}%
m_{\lambda},
\]
for $j=0,1,2,\ldots$.
\end{definition}

We establish some basic properties as well as an explicit formula for this
coefficient, and also explain the relation with the generalized binomial coefficient.

\begin{proposition}
\label{bincof}For $\lambda,\nu\in\mathcal{P}^{\left(  N\right)  }$ the
following hold:

\begin{enumerate}
\item  let $j=\left|  \lambda\right|  -\left|  \nu\right|  >0$, then
\[%
\genfrac{\langle}{\rangle}{0pt}{}{\lambda}{\nu}%
=\sum\left\{  \binom{j}{\lambda-\sigma}:\sigma^{+}=\nu\right\}  ;
\]
note $\binom{j}{\lambda-\sigma}=0$ if any component $\lambda_{i}-\sigma_{i}%
<0$, the sum is over all distinct permutations of $\nu$, (including any
possible zero components to make $\nu$ an element of $\mathbb{N}_{0}^{N}$);

\item $%
\genfrac{\langle}{\rangle}{0pt}{}{\lambda}{\nu}%
\neq0$ if and only if $\nu\subset\lambda$;

\item $%
\genfrac{\langle}{\rangle}{0pt}{}{\lambda}{\lambda}%
=1$.
\end{enumerate}
\end{proposition}

\begin{proof}
Expanding $e_{1}^{j}m_{\nu}$ we obtain
\[
e_{1}^{j}\,m_{\nu}=\sum\left\{  \binom{j}{\alpha}x^{\alpha}x^{\sigma}%
:\alpha\in\mathbb{N}_{0}^{N},\left|  \alpha\right|  =j,\sigma^{+}=\nu\right\}
.
\]
Since this is a symmetric polynomial, the coefficient of $m_{\lambda}$ equals
the coefficient of $x^{\lambda}$, namely the sum of $\binom{j}{\alpha}$ with
$\alpha+\sigma=\lambda$ (component-wise addition). For any $\sigma$ with
$\sigma^{+}=\nu$ if $\sigma_{i}>\lambda_{i}$ for some $i$ then the
corresponding term is zero, by definition of the multinomial coefficient. For
part (2), $\nu\subset\lambda$ (and $\left|  \lambda\right|  =\left|
\nu\right|  +j$) implies $\binom{j}{\lambda-\nu}>0$; conversely if $\binom
{j}{\lambda-\sigma}\neq0$ for some $\sigma$ with $\sigma^{+}=\nu$ then
$\nu\subset\lambda$ (if $\sigma_{i}\leq\lambda_{i}$ for each $i$ then
$\sigma^{+}\subset\lambda$; indeed, let $w$ be a permutation so that $\nu
_{i}=\sigma_{w\left(  i\right)  }$ for each $i$, so that $\nu_{j}\leq\nu
_{i}=\sigma_{w\left(  i\right)  }\leq\lambda_{w\left(  i\right)  }$ for $1\leq
i\leq j$ then $\nu_{j}\leq\lambda_{j}$ because it is less than or equal to at
least $j$ components of the partition $\lambda$). Part (3) is trivial.
\end{proof}

The term ``generalized binomial coefficient'' has appeared in several contexts
(such as Jack polynomials, see Lassalle \cite{L}; part 2 of Proposition
\ref{bincof} is a special case of his results); here we use the version
defined by
\begin{equation}
\frac{m_{\lambda}\left(  x+s1^{N}\right)  }{m_{\lambda}\left(  1^{N}\right)
}=\sum_{\nu}s^{\left|  \lambda\right|  -\left|  \nu\right|  }\binom{\lambda
}{\nu}\frac{m_{\nu}\left(  x\right)  }{m_{\nu}\left(  1^{N}\right)  }%
,s\in\mathbb{R}. \label{genbin}%
\end{equation}
This is essentially the Taylor series for $m_{\lambda}$ evaluated at $\left(
x_{1}+s,\ldots,x_{N}+s\right)  $. The coefficient $\dbinom{\lambda}{\nu
}=\dfrac{\lambda!}{\nu!\left(  \left|  \lambda\right|  -\left|  \nu\right|
\right)  !}%
\genfrac{\langle}{\rangle}{0pt}{0}{\lambda}{\nu}%
$. This formula is established by means of the inner product $\left\langle
f,g\right\rangle =f\left(  \partial\right)  g\left(  x\right)  |_{x=0}$; apply
$\left(  \frac{\partial}{\partial s}\right)  ^{j}$ to both sides of equation
\ref{genbin} and set $s=0$ to show
\[
\frac{1}{m_{\lambda}\left(  1^{N}\right)  }\left(  \sum_{i=1}^{N}%
\frac{\partial}{\partial x_{i}}\right)  ^{j}m_{\lambda}\left(  x\right)
=\sum_{\left|  \nu\right|  =\left|  \lambda\right|  -j}j!\binom{\lambda}{\nu
}\frac{m_{\nu}\left(  x\right)  }{m_{\nu}\left(  1^{N}\right)  };
\]
take the inner product with $m_{\sigma}$ for any $\sigma\in\mathcal{P}%
^{\left(  N\right)  }$ with $\left|  \sigma\right|  =\left|  \lambda\right|
-j$; note that $\left\langle m_{\sigma},m_{\nu}\right\rangle =\delta
_{\sigma\nu}\sigma!m_{\sigma}\left(  1^{N}\right)  $ and $\left\langle \left(
\sum_{i=1}^{N}\frac{\partial}{\partial x_{i}}\right)  ^{j}m_{\lambda
},m_{\sigma}\right\rangle =\left\langle m_{\lambda},e_{1}^{j}m_{\sigma
}\right\rangle $. The coefficient $\binom{\lambda}{\nu}$ is thus shown to be
independent of $N$, provided the number of (nonzero) parts of $\lambda,\nu$
does not exceed $N$. There are useful formulae involving contiguous partitions
(if $\nu\supset\lambda$ and $\left|  \nu\right|  =\left|  \lambda\right|  +1$
then $\nu=\lambda+\varepsilon_{i}$ for $i=1$ or $\lambda_{i-1}>\lambda_{i}$).

\begin{proposition}
\label{bin1}Suppose $\lambda,\sigma\in\mathcal{P}$ and $\left|  \sigma\right|
\geq\left|  \lambda\right|  +1$, also $j=1$ or $\lambda_{j-1}>\lambda_{j}$
then
\begin{align}%
\genfrac{\langle}{\rangle}{0pt}{}{\lambda+\varepsilon_{j}}{\lambda}%
&  =1+\#\left\{  i:\lambda_{i}=\lambda_{j}+1\right\}  ,\label{bin1a}\\%
\genfrac{\langle}{\rangle}{0pt}{}{\sigma}{\lambda}%
&  =\sum\left\{
\genfrac{\langle}{\rangle}{0pt}{}{\sigma}{\lambda+\varepsilon_{i}}%
\genfrac{\langle}{\rangle}{0pt}{}{\lambda+\varepsilon_{i}}{\lambda}%
:i=1\text{ or }\lambda_{i-1}>\lambda_{i}\right\}  . \label{bin1b}%
\end{align}
\end{proposition}

\begin{proof}
Suppose $\lambda_{s-1}>\lambda_{s}=\lambda_{j-1}=\lambda_{j}+1$ (or $s=1, $or
$j=1)$, then there are exactly $j-s+1$ distinct permutations $w\lambda$ for
which $\lambda+\varepsilon_{j}-w\lambda$ has no negative components (namely
the transpositions $w=\left(  i,j\right)  $ for some $i$ with $s\leq i\leq
j-1$, or $w=1$). Each such term contributes $1$ to the sum defining $%
\genfrac{\langle}{\rangle}{0pt}{}{\lambda+\varepsilon_{j}}{\lambda}%
$. For the second part, let $n=$ $\left|  \sigma\right|  -\left|
\lambda\right|  $ and consider the coefficient of $m_{\sigma}$ in
\begin{align*}
e_{1}^{n}m_{\lambda}  &  =\sum_{\left|  \nu\right|  =\left|  \lambda\right|
+n}%
\genfrac{\langle}{\rangle}{0pt}{}{\nu}{\lambda}%
m_{\nu}=e_{1}^{n-1}\sum_{i}%
\genfrac{\langle}{\rangle}{0pt}{}{\lambda+\varepsilon_{i}}{\lambda}%
m_{\lambda+\varepsilon_{i}}\\
&  =\sum_{\left|  \nu\right|  =\left|  \lambda\right|  +n}\sum_{i}%
\genfrac{\langle}{\rangle}{0pt}{}{\nu}{\lambda+\varepsilon_{i}}%
\genfrac{\langle}{\rangle}{0pt}{}{\lambda+\varepsilon_{i}}{\lambda}%
m_{\nu},
\end{align*}
summing over $i=1$ and $i\geq2$ with $\lambda_{i-1}>\lambda_{i}.$
\end{proof}

We are ready to consider the transition matrices between the $m_{\lambda}$ and
$\widetilde{m}_{\lambda}$ bases (that the second set is a basis will be
immediately proven).

\begin{definition}
The transition matrices $A\left(  \mu,\lambda\right)  ,B\left(  \mu
,\lambda\right)  $ (for $\lambda,\mu\in\mathcal{P},\,\left|  \lambda\right|
=\left|  \mu\right|  $) are defined by
\begin{align*}
\widetilde{m}_{\lambda}  &  =\sum_{\mu}A\left(  \mu,\lambda\right)  m_{\mu},\\
m_{\lambda}  &  =\sum_{\mu}B\left(  \mu,\lambda\right)  \widetilde{m}_{\mu}.
\end{align*}
\end{definition}

By Definition \ref{e1m} it follows that
\[
A\left(  \mu,\lambda\right)  =%
\genfrac{\langle}{\rangle}{0pt}{}{\mu}{\left(  \lambda_{2},\lambda_{2}%
,\lambda_{3},\ldots\right)  }%
.
\]
The condition $A\left(  \mu,\lambda\right)  \neq0$ leads to the following
ordering of partitions:

\begin{definition}
For $\lambda,\mu\in\mathcal{P}$ the relation $\mu\preceq\lambda$ means that
$\left|  \lambda\right|  =\left|  \mu\right|  $ and $\lambda_{i}\leq\mu_{i}$
for all $i\geq2$.
\end{definition}

It\ is clear that $\preceq$ is a partial ordering, $\mu\preceq\lambda$ implies
$\mu_{1}\leq\lambda_{1},l\left(  \lambda\right)  \leq l\left(  \mu\right)  $
and that the maximum element of $\mathcal{P}_{n}^{\left(  N\right)  } $ is
$\left(  n\right)  $. In fact $\mu\preceq\lambda$ implies that $\lambda$
dominates $\mu$, but $\preceq$ is not identical to the dominance ordering.
Further the partitions $\left\{  \lambda:\lambda_{1}=\lambda_{2}\right\}  $
are minimal elements; indeed suppose $\mu\preceq\lambda$ and $\lambda
_{1}=\lambda_{2}$, then $\mu_{1}\leq\lambda_{1}=\lambda_{2}\leq\mu_{2}$ which
implies $\mu_{1}=\mu_{2}=\lambda_{1}$; thus $\sum_{i\geq2}\left(  \mu
_{i}-\lambda_{i}\right)  =0$ and $\lambda=\mu$ (by definition, $\mu
_{i}-\lambda_{i}\geq0$ for each $i\geq2$). There are exceptional minimal
elements of the form $\left(  m+1,m,\ldots,m\right)  $, with $n=Nm+1$, but
these are minimal only for $\mathcal{P}^{\left(  N\right)  } $ since $\left(
m,m,\ldots,m,1\right)  \preceq\left(  m+1,m,\ldots,m,0\right)  $. We see that
the matrix $A\left(  \mu,\lambda\right)  $ is triangular ($A\left(
\mu,\lambda\right)  \neq0$ implies $\mu\preceq\lambda$), unipotent ($A\left(
\lambda,\lambda\right)  =1$) and its entries are nonnegative integers.

\begin{proposition}
For $n\geq0$ the set $\left\{  \widetilde{m}_{\lambda}:\lambda\in
\mathcal{P}_{n}^{\left(  N\right)  }\right\}  $ is a basis for the symmetric
polynomials of degree $n$ in $N$ variables. The transition matrix $B\left(
\mu,\lambda\right)  $ is triangular for the ordering $\preceq$ and is
unipotent with integer entries.
\end{proposition}

Regarding the dependence on $N$: begin by assuming that $n\leq N$ so that
$l\left(  \lambda\right)  \leq N$ for each $\lambda\in\mathcal{P}_{n}$ and the
defining equations for $A\left(  \mu,\lambda\right)  ,B\left(  \mu
,\lambda\right)  $ are unambiguous, then to restrict to a smaller number of
variables, say $M<N$, substitute $x_{M+1}=x_{M+2}=\ldots=x_{N}=0$ with the
effect of removing all terms with $l\left(  \lambda\right)  >M$ or $l\left(
\mu\right)  >M$. Thus the transition matrices for the case of $M$ variables
with $n>M$ are principal submatrices of the general transition matrices
(deleting rows and columns labeled by $\lambda$ with $l\left(  \lambda\right)
>M$). If the index set $\mathcal{P}_{n}^{\left(  N\right)  }$ is ordered first
in decreasing order of $l\left(  \lambda\right)  $ and then with respect to
$\preceq$ (possible, since $\mu\preceq\lambda$ implies $l\left(  \mu\right)
\geq l\left(  \lambda\right)  $) then the triangularity makes it obvious that
the submatrix of $B$ is the inverse of the submatrix of $A$; see the example
of $B$ for $\mathcal{P}_{6}^{\left(  4\right)  }$ at the end of this section.

It remains to establish an explicit formula for $B\left(  \mu,\lambda\right)
$. The formula itself was postulated based on computer algebra
experimentation, and explicit known formulae for $N=3$. The proof will be by
induction on $n$ and requires showing that the claimed formula satisfies the
following contiguity relations.

\begin{proposition}
\label{bcontig}For $\lambda,\mu\in\mathcal{P}_{n}^{\left(  N\right)  }$, and
$\sigma\in\mathcal{P}_{n+1}^{\left(  N\right)  }$ with $\sigma_{1}=\sigma_{2}%
$,
\begin{align}
\sum_{\nu\supset\lambda,\left|  \nu\right|  =\left|  \lambda\right|  +1}%
\genfrac{\langle}{\rangle}{0pt}{}{\nu}{\lambda}%
B\left(  \mu+\varepsilon_{1},\nu\right)   &  =B\left(  \mu,\lambda\right)
,\label{bform}\\
\sum_{\nu\supset\lambda,\left|  \nu\right|  =\left|  \lambda\right|  +1}%
\genfrac{\langle}{\rangle}{0pt}{}{\nu}{\lambda}%
B\left(  \sigma,\nu\right)   &  =0. \label{bform2}%
\end{align}
\end{proposition}

\begin{proof}
Multiply both sides of the equation $m_{\lambda}=\sum_{\mu\preceq\lambda
}B\left(  \mu,\lambda\right)  \widetilde{m}_{\mu}$ by $e_{1}$ to obtain
\begin{align*}
\sum_{\nu\supset\lambda,\left|  \nu\right|  =\left|  \lambda\right|  +1}%
\genfrac{\langle}{\rangle}{0pt}{}{\nu}{\lambda}%
m_{\nu}  &  =\sum_{\mu\preceq\lambda}B\left(  \mu,\lambda\right)
\widetilde{m}_{\mu+\varepsilon_{1}}\\
&  =\sum_{\nu\supset\lambda,\left|  \nu\right|  =\left|  \lambda\right|  +1}%
\genfrac{\langle}{\rangle}{0pt}{}{\nu}{\lambda}%
\sum_{\tau\preceq\nu}B\left(  \tau,\nu\right)  \widetilde{m}_{\tau}.
\end{align*}
Since the set $\left\{  \widetilde{m}_{\tau}:\tau\in\mathcal{P}_{n+1}^{\left(
N\right)  }\right\}  $ is a basis the desired equations are consequences of
matching coefficients in the two right-hand expansions.
\end{proof}

\begin{lemma}
Suppose $B^{\prime}\left(  \mu,\lambda\right)  $ is a matrix satisfying the
equations \ref{bform},\ref{bform2} with $B$ replaced by $B^{\prime}$, has the
same triangularity property as $B,$and $B^{\prime}\left(  \sigma
,\sigma\right)  =1$ for each $\sigma\in\mathcal{P}$ with $\sigma_{1}%
=\sigma_{2}$, then $B^{\prime}\left(  \mu,\lambda\right)  =B\left(
\mu,\lambda\right)  $ for all $\mu,\lambda\in\mathcal{P}$ (with $\left|
\lambda\right|  =\left|  \mu\right|  $).
\end{lemma}

\begin{proof}
This is a double induction on $\left|  \lambda\right|  $ and $\lambda
_{1}-\lambda_{2}.$ For any $\sigma\in\mathcal{P}$ with $\sigma_{1}=\sigma_{2}$
the hypothesis shows that $B^{\prime}\left(  \mu,\sigma\right)  =0$ for
$\mu\neq\sigma$ (since $\sigma$ is minimal) and $B^{\prime}\left(
\sigma,\sigma\right)  =1$, thus $B^{\prime}\left(  \mu,\sigma\right)
=B\left(  \mu,\sigma\right)  $ for all $\mu$ (with $\left|  \mu\right|
=\left|  \sigma\right|  $). Suppose that $B^{\prime}\left(  \mu,\lambda
\right)  =B\left(  \mu,\lambda\right)  $ for all $\mu,\lambda\in\mathcal{P}$
with $\left|  \lambda\right|  =\left|  \mu\right|  =n$ or with $\left|
\lambda\right|  =\left|  \mu\right|  =n+1$ and $\lambda_{1}-\lambda_{2}\leq j$
for some $j\geq0$. Fix $\nu\in\mathcal{P}$ with $\left|  \nu\right|  =n+1$ and
$\nu_{1}-\nu_{2}=j+1$ and let $\lambda=\nu-\varepsilon_{1}$. Then $\tau
\supset\lambda$ and $\left|  \tau\right|  =n+1$ implies $\tau=\nu$ or
$\tau=\lambda+\varepsilon_{i}$ with $i\geq2$ and $\lambda_{i-1}>\lambda
_{i}=\nu_{i}$. Note that $%
\genfrac{\langle}{\rangle}{0pt}{}{\lambda+\varepsilon_{1}}{\lambda}%
=1$ by equation \ref{bin1a}. For any $\mu\in\mathcal{P}_{n+1}$, by hypothesis
\[
B^{\prime}\left(  \mu,\nu\right)  =B^{\prime}\left(  \mu-\varepsilon
_{1},\lambda\right)  -\sum_{i\geq2,\lambda_{i-1}>\lambda_{i}}%
\genfrac{\langle}{\rangle}{0pt}{}{\lambda+\varepsilon_{i}}{\lambda}%
B^{\prime}\left(  \mu,\lambda+\varepsilon_{i}\right)  ,
\]
replacing $B^{\prime}\left(  \mu-\varepsilon_{1},\lambda\right)  $ by $0$ if
$\mu_{1}=\mu_{2}$. By the inductive hypothesis each term $B^{\prime}\left(
\sigma,\tau\right)  $ appearing on the right hand side satisfies $B^{\prime
}\left(  \sigma,\tau\right)  =B\left(  \sigma,\tau\right)  $. By the
Proposition, $B^{\prime}\left(  \mu,\nu\right)  =B\left(  \mu,\nu\right)  $.
This completes the induction.
\end{proof}

In the following formula the coefficient $%
\genfrac{\langle}{\rangle}{0pt}{}{\sigma}{\tau}%
$ is used with partitions whose first parts are deleted, also $\left(  \mu
_{2}-1,\mu_{3},\ldots\right)  $ is not a partition if $\mu_{2}=\mu_{3}$ and so
the $\alpha^{+}$ notation is used; the binomial coefficient $\binom{-1}{m}=0$.

\begin{theorem}
\label{bformula}For $\mu,\lambda\in\mathcal{P}_{n}$%
\begin{align*}
B\left(  \mu,\lambda\right)   &  =\left(  -1\right)  ^{\lambda_{1}-\mu_{1}%
}\times\\
&  \left\{  \binom{\lambda_{1}-\mu_{2}}{\mu_{1}-\mu_{2}}%
\genfrac{\langle}{\rangle}{0pt}{}{\left(  \mu_{2},\mu_{3},\ldots\right)
}{\left(  \lambda_{2},\lambda_{3},\ldots\right)  }%
+\binom{\lambda_{1}-\mu_{2}-1}{\mu_{1}-\mu_{2}}%
\genfrac{\langle}{\rangle}{0pt}{}{\left(  \mu_{2}-1,\mu_{3},\ldots\right)
^{+}}{\left(  \lambda_{2},\lambda_{3},\ldots\right)  }%
\right\}
\end{align*}
\end{theorem}

\begin{proof}
The proof is broken down into cases depending on the values of $\lambda
_{1}-\mu_{1}$ and $\mu_{1}-\mu_{2}$. We let $B^{\prime}\left(  \mu
,\lambda\right)  $ denote the right-hand-side of the stated formula and
proceed as in the Lemma to show $B^{\prime}$ satisfies equations \ref{bform}
and \ref{bform2}. It is clear that $B^{\prime}\left(  \mu,\lambda\right)
\neq0$ implies $\lambda_{i}\leq\mu_{i}$ for all $i\geq2$, that is, $\mu
\preceq\lambda$. Also $B^{\prime}\left(  \mu,\mu\right)  =1$ for each $\mu
\in\mathcal{P}_{n}$. Fix $\mu,\lambda\in\mathcal{P}_{n}$ and let $\mu^{\prime
}=\left(  \mu_{2},\mu_{3},\ldots\right)  ,\mu^{\prime\prime}=\left(  \mu
_{2}-1,\mu_{3},\ldots\right)  ^{+},\allowbreak\sigma^{\prime}=\left(
\sigma_{2},\sigma_{3},\ldots\right)  ,\allowbreak\sigma^{\prime\prime}=\left(
\sigma_{2}-1,\sigma_{3},\ldots\right)  ^{+},\allowbreak\lambda^{\prime
}=\left(  \lambda_{2},\lambda_{3},\ldots\right)  $ (these are partitions, but
with components labeled by $i\geq2$; for example $\lambda^{\prime}%
+\varepsilon_{2}=\left(  \lambda_{2}+1,\lambda_{3},\ldots\right)  $). The
partitions $\nu$ satisfying $\nu\supset\lambda$ and $\left|  \nu\right|
=\left|  \lambda\right|  +1$ are of the form $\lambda+\varepsilon_{i}$ with
$i=1$ or $\lambda_{i-1}>\lambda_{i}$. We rewrite equation \ref{bform} (not yet
proven!) as
\[
\sum_{i\geq2,\lambda_{i-1}>\lambda_{i}}%
\genfrac{\langle}{\rangle}{0pt}{}{\lambda+\varepsilon_{i}}{\lambda}%
B^{\prime}\left(  \mu+\varepsilon_{1},\lambda+\varepsilon_{i}\right)
=B^{\prime}\left(  \mu,\lambda\right)  -B^{\prime}\left(  \mu+\varepsilon
_{1},\lambda+\varepsilon_{1}\right)  .
\]
The left-hand-side equals
\begin{align*}
&  \left(  -1\right)  ^{\lambda_{1}-\mu_{1}-1}\times\\
&  \sum_{i\geq2,\lambda_{i-1}>\lambda_{i}}\left\{  \binom{\lambda_{1}-\mu_{2}%
}{\mu_{1}+1-\mu_{2}}%
\genfrac{\langle}{\rangle}{0pt}{}{\mu^{\prime}}{\lambda^{\prime}%
+\varepsilon_{i}}%
+\binom{\lambda_{1}-\mu_{2}-1}{\mu_{1}+1-\mu_{2}}%
\genfrac{\langle}{\rangle}{0pt}{}{\mu^{\prime\prime}}{\lambda^{\prime
}+\varepsilon_{i}}%
\right\}
\genfrac{\langle}{\rangle}{0pt}{}{\lambda+\varepsilon_{i}}{\lambda}%
.
\end{align*}
The right-hand-side equals
\begin{align*}
&  \left(  -1\right)  ^{\lambda_{1}-\mu_{1}}\left\{  \left(  \binom
{\lambda_{1}-\mu_{2}}{\mu_{1}-\mu_{2}}-\binom{\lambda_{1}+1-\mu_{2}}{\mu
_{1}+1-\mu_{2}}\right)
\genfrac{\langle}{\rangle}{0pt}{}{\mu^{\prime}}{\lambda^{\prime}}%
\right. \\
&  +\left.  \left(  \binom{\lambda_{1}-\mu_{2}-1}{\mu_{1}-\mu_{2}}%
-\binom{\lambda_{1}-\mu_{2}}{\mu_{1}+1-\mu_{2}}\right)
\genfrac{\langle}{\rangle}{0pt}{}{\mu^{\prime\prime}}{\lambda^{\prime}}%
\right\} \\
&  =\left(  -1\right)  ^{\lambda_{1}-\mu_{1}-1}\left(  \binom{\lambda_{1}%
-\mu_{2}}{\mu_{1}+1-\mu_{2}}%
\genfrac{\langle}{\rangle}{0pt}{}{\mu^{\prime}}{\lambda^{\prime}}%
+\binom{\lambda_{1}-\mu_{2}-1}{\mu_{1}+1-\mu_{2}}%
\genfrac{\langle}{\rangle}{0pt}{}{\mu^{\prime\prime}}{\lambda^{\prime}}%
\right)  .
\end{align*}
Similarly rewrite equation \ref{bform2} as
\[
\sum_{i\geq2,\lambda_{i-1}>\lambda_{i}}%
\genfrac{\langle}{\rangle}{0pt}{}{\lambda+\varepsilon_{i}}{\lambda}%
B^{\prime}\left(  \sigma,\lambda+\varepsilon_{i}\right)  =-B^{\prime}\left(
\sigma,\lambda+\varepsilon_{1}\right)  ,
\]
with the left-hand-side (of course, $\sigma_{1}=\sigma_{2}$)
\begin{align*}
&  \left(  -1\right)  ^{\lambda_{1}-\sigma_{1}}\\
&  \times\sum_{i\geq2,\lambda_{i-1}>\lambda_{i}}\left\{  \binom{\lambda
_{1}-\sigma_{2}}{\sigma_{1}-\sigma_{2}}%
\genfrac{\langle}{\rangle}{0pt}{}{\sigma^{\prime}}{\lambda^{\prime
}+\varepsilon_{i}}%
+\binom{\lambda_{1}-\sigma_{2}-1}{\sigma_{1}-\sigma_{2}}%
\genfrac{\langle}{\rangle}{0pt}{}{\sigma^{\prime\prime}}{\lambda^{\prime
}+\varepsilon_{i}}%
\right\}
\genfrac{\langle}{\rangle}{0pt}{}{\lambda+\varepsilon_{i}}{\lambda}%
,
\end{align*}
and the right-hand-side
\[
\left(  -1\right)  ^{\lambda_{1}-\sigma_{1}}\left(  \binom{\lambda_{1}%
-\sigma_{2}}{\sigma_{1}-\sigma_{2}}%
\genfrac{\langle}{\rangle}{0pt}{}{\sigma^{\prime}}{\lambda^{\prime}}%
+\binom{\lambda_{1}-\sigma_{2}-1}{\sigma_{1}-\sigma_{2}}%
\genfrac{\langle}{\rangle}{0pt}{}{\sigma^{\prime\prime}}{\lambda^{\prime}}%
\right)  .
\]
The goal is to show that the two sides are equal. We reduce the two cases to
one by replacing $\sigma_{1},\sigma^{\prime},\sigma^{\prime\prime}$ by
$\mu_{1}+1,\mu^{\prime},\mu^{\prime\prime}$respectively. Let $E$ denote the
difference between the two sides (left - right), and let
\begin{align*}
E_{0}  &  =-\binom{\lambda_{1}-\mu_{2}}{\mu_{1}+1-\mu_{2}}%
\genfrac{\langle}{\rangle}{0pt}{}{\mu^{\prime}}{\lambda^{\prime}}%
,\\
E_{1}  &  =\binom{\lambda_{1}-\mu_{2}}{\mu_{1}+1-\mu_{2}}\sum_{i\geq
2,\lambda_{i-1}>\lambda_{i}}%
\genfrac{\langle}{\rangle}{0pt}{}{\mu^{\prime}}{\lambda^{\prime}%
+\varepsilon_{i}}%
\genfrac{\langle}{\rangle}{0pt}{}{\lambda+\varepsilon_{i}}{\lambda}%
-\binom{\lambda_{1}-\mu_{2}-1}{\mu_{1}+1-\mu_{2}}%
\genfrac{\langle}{\rangle}{0pt}{}{\mu^{\prime\prime}}{\lambda^{\prime}}%
,\\
E_{2}  &  =\binom{\lambda_{1}-\mu_{2}-1}{\mu_{1}+1-\mu_{2}}\sum_{i\geq
2,\lambda_{i-1}>\lambda_{i}}%
\genfrac{\langle}{\rangle}{0pt}{}{\mu^{\prime\prime}}{\lambda^{\prime
}+\varepsilon_{i}}%
\genfrac{\langle}{\rangle}{0pt}{}{\lambda+\varepsilon_{i}}{\lambda}%
,
\end{align*}
then $E=\left(  -1\right)  ^{\lambda_{1}-\mu_{1}-1}\left(  E_{0}+E_{1}%
+E_{2}\right)  $. We can assume $\mu^{\prime}\supset\lambda^{\prime}$
($\mu\preceq\lambda$), because if any term in $E$ is nonzero then $%
\genfrac{\langle}{\rangle}{0pt}{}{\mu^{\prime}}{\lambda^{\prime}}%
\neq0$ (if $%
\genfrac{\langle}{\rangle}{0pt}{}{\mu^{\prime\prime}}{\lambda^{\prime
}+\varepsilon_{i}}%
\neq0$ for some $i$ then $\mu^{\prime}\supset\mu^{\prime\prime}\supset
\lambda^{\prime}+\varepsilon_{i}\supset\lambda^{\prime}$, and similarly if $%
\genfrac{\langle}{\rangle}{0pt}{}{\mu^{\prime}}{\lambda^{\prime}%
+\varepsilon_{i}}%
\neq0$ for some $i$ then $\mu^{\prime}\supset\lambda^{\prime}$).

\begin{enumerate}
\item \textbf{Case}: $\left|  \mu^{\prime}\right|  =\left|  \lambda^{\prime
}\right|  $: Here $E_{1}=0=E_{2}$, because any coefficient $%
\genfrac{\langle}{\rangle}{0pt}{}{\sigma}{\tau}%
$ with $\left|  \sigma\right|  <\left|  \tau\right|  $ is zero. Further
$\mu_{1}=n-\left|  \mu^{\prime}\right|  =\lambda_{1}$ and thus $E_{0}=0$ (from
the binomial coefficient).

\item \textbf{Case}: $\left|  \mu^{\prime}\right|  =\left|  \lambda^{\prime
}\right|  +1$ and $\mu_{1}\geq\mu_{2}$: Here $E_{2}=0$ and $\lambda_{1}%
=\mu_{1}+1,$thus the coefficient of $%
\genfrac{\langle}{\rangle}{0pt}{}{\mu^{\prime\prime}}{\lambda^{\prime}}%
$ is zero (specifically, $\binom{\lambda_{1}-\mu_{2}-1}{\mu_{1}-\mu_{2}%
}-\binom{\lambda_{1}-\mu_{2}}{\mu_{1}+1-\mu_{2}}=1-1$). The hypothesis implies
$\mu^{\prime}=\lambda^{\prime}+\varepsilon_{j}$ for some $j\geq2$, that is,
$\mu=\left(  \lambda_{1}-1,\lambda_{2},\ldots,\lambda_{j}+1,\ldots\right)  $.
If $j=2$ then $\lambda_{1}\geq\lambda_{2}+2$ otherwise $j>2$ and $\lambda
_{1}-1\geq\lambda_{j-1}>\lambda_{j}$. Thus $E_{0}+E_{1}=%
\genfrac{\langle}{\rangle}{0pt}{}{\lambda+\varepsilon_{j}}{\lambda}%
-%
\genfrac{\langle}{\rangle}{0pt}{}{\lambda^{\prime}+\varepsilon_{j}%
}{\lambda^{\prime}}%
$; only the term with $i=j$ in the sum is nonzero (we showed $\lambda
_{j-1}>\lambda_{j}$). By equation \ref{bin1a} $%
\genfrac{\langle}{\rangle}{0pt}{}{\lambda+\varepsilon_{j}}{\lambda}%
=%
\genfrac{\langle}{\rangle}{0pt}{}{\lambda^{\prime}+\varepsilon_{j}%
}{\lambda^{\prime}}%
$, since $1\notin\left\{  i:\lambda_{i}=\lambda_{j}+1\right\}  $.

\item \textbf{Case}: $\left|  \mu^{\prime}\right|  =\left|  \lambda^{\prime
}\right|  +1$ and $\mu_{1}=\mu_{2}-1$ (that is, $\mu=\sigma-\varepsilon_{1}$
where $\sigma_{1}=\sigma_{2}$): Here $E_{2}=0$ and the binomial coefficients
in $E_{0},E_{1}$ are all equal to $1$. Let $\mu^{\prime}=\lambda^{\prime
}+\varepsilon_{j}$. If $\mu^{\prime\prime}=\lambda^{\prime}$ then either
$j=2,\,\lambda=\left(  \lambda_{1},\lambda_{1}-1,\lambda_{3},\ldots\right)
,\,\mu=\left(  \lambda_{1}-1,\lambda_{1},\lambda_{3},\ldots\right)  $ with
$\lambda_{1}>\lambda_{3}$ or $j>2,\,\lambda=\left(  \lambda_{1},\lambda
_{1},\ldots,\lambda_{1},\lambda_{1}-1,\lambda_{j+1},\ldots\right)
,\allowbreak\,\mu=\left(  \lambda_{1}-1,\lambda_{1},\ldots,\lambda_{1}%
,\lambda_{j+1},\ldots\right)  $ (that is, $\lambda_{i}=\lambda_{1}$ for $1\leq
i\leq j-1$ and $\lambda_{j}=\lambda_{1}-1$), which implies
\begin{align*}
E_{0}+E_{1}  &  =%
\genfrac{\langle}{\rangle}{0pt}{0}{\mu^{\prime}}{\lambda^{\prime}%
+\varepsilon_{j}}%
\genfrac{\langle}{\rangle}{0pt}{0}{\lambda+\varepsilon_{j}}{\lambda}%
-%
\genfrac{\langle}{\rangle}{0pt}{0}{\mu^{\prime\prime}}{\lambda^{\prime}}%
-%
\genfrac{\langle}{\rangle}{0pt}{0}{\mu^{\prime}}{\lambda^{\prime}}%
\\
&  =%
\genfrac{\langle}{\rangle}{0pt}{0}{\lambda+\varepsilon_{j}}{\lambda}%
-1-%
\genfrac{\langle}{\rangle}{0pt}{0}{\lambda^{\prime}+\varepsilon_{j}%
}{\lambda^{\prime}}%
=j-1-(j-1),
\end{align*}
by equation \ref{bin1a}. If $\mu^{\prime\prime}\neq\lambda^{\prime}$ then
$\lambda_{j}+1<\lambda_{2}$ (so $j>2)$ and $E_{0}+E_{1}=%
\genfrac{\langle}{\rangle}{0pt}{}{\lambda+\varepsilon_{j}}{\lambda}%
-%
\genfrac{\langle}{\rangle}{0pt}{}{\lambda^{\prime}+\varepsilon_{j}%
}{\lambda^{\prime}}%
=0$.

\item \textbf{Case}: $\left|  \mu^{\prime}\right|  \geq\left|  \lambda
^{\prime}\right|  +2$ and $\mu_{1}\geq\mu_{2}$: Thus $\lambda_{1}\geq\mu
_{1}+2\geq\mu_{2}+2\geq\lambda_{2}+2$ and equation \ref{bin1a} implies $%
\genfrac{\langle}{\rangle}{0pt}{}{\lambda+\varepsilon_{i}}{\lambda}%
=%
\genfrac{\langle}{\rangle}{0pt}{}{\lambda^{\prime}+\varepsilon_{i}%
}{\lambda^{\prime}}%
$ for each $i\geq2,\lambda_{i-1}>\lambda_{i}$ (which includes $i=2$). Also
equation \ref{bin1b} applied to the truncated partition $\tau^{\prime}=\left(
\tau_{2},\tau_{3},\ldots\right)  $ with $\left|  \tau^{\prime}\right|
>\left|  \lambda^{\prime}\right|  $ shows that
\[
\sum\left\{
\genfrac{\langle}{\rangle}{0pt}{}{\tau^{\prime}}{\lambda^{\prime}%
+\varepsilon_{i}}%
\genfrac{\langle}{\rangle}{0pt}{}{\lambda^{\prime}+\varepsilon_{i}%
}{\lambda^{\prime}}%
:i=2\text{ or }i>2,\lambda_{i-1}>\lambda\right\}  =%
\genfrac{\langle}{\rangle}{0pt}{}{\tau^{\prime}}{\lambda^{\prime}}%
.
\]
Substituting $\tau^{\prime}=\mu^{\prime}$ and $\tau^{\prime}=\mu^{\prime
\prime}$ in this equation shows that $E_{0}+E_{1}+E_{2}=0.$

\item \textbf{Case}: $\left|  \mu^{\prime}\right|  \geq\left|  \lambda
^{\prime}\right|  +2$ and $\mu_{1}=\mu_{2}-1$ (that is, $\mu=\sigma
-\varepsilon_{1}$ where $\sigma_{1}=\sigma_{2}$): Similarly $\lambda_{1}%
\geq\left(  \mu_{2}-1\right)  +2\geq\lambda_{2}+1$, and $%
\genfrac{\langle}{\rangle}{0pt}{}{\lambda+\varepsilon_{i}}{\lambda}%
=%
\genfrac{\langle}{\rangle}{0pt}{}{\lambda^{\prime}+\varepsilon_{i}%
}{\lambda^{\prime}}%
$ for each $i>2$ with $\lambda_{i-1}>\lambda_{i}$.\ Either $\lambda_{1}%
\geq\lambda_{2}+2$ which implies $%
\genfrac{\langle}{\rangle}{0pt}{}{\lambda+\varepsilon_{2}}{\lambda}%
=%
\genfrac{\langle}{\rangle}{0pt}{}{\lambda^{\prime}+\varepsilon_{2}%
}{\lambda^{\prime}}%
=1$ or $\lambda_{1}=\lambda_{2}+1$ which implies $\mu_{2}=\lambda_{2}$ and $%
\genfrac{\langle}{\rangle}{0pt}{}{\mu^{\prime}}{\lambda^{\prime}%
+\varepsilon_{2}}%
=0=%
\genfrac{\langle}{\rangle}{0pt}{}{\mu^{\prime\prime}}{\lambda^{\prime
}+\varepsilon_{2}}%
$ because $\left(  \lambda^{\prime}+\varepsilon_{2}\right)  _{2}>\mu_{2}%
\geq\left(  \mu^{\prime\prime}\right)  _{2}$. In both cases equation
\ref{bin1b} shows $E_{0}+E_{1}+E_{2}=0$ as previously (in the case
$\lambda_{1}=\lambda_{2}+1$ one has $%
\genfrac{\langle}{\rangle}{0pt}{}{\lambda+\varepsilon_{2}}{\lambda}%
=2,%
\genfrac{\langle}{\rangle}{0pt}{}{\lambda^{\prime}+\varepsilon_{2}%
}{\lambda^{\prime}}%
=1$).
\end{enumerate}

Thus the matrix $B^{\prime}\left(  \lambda,\mu\right)  $ satisfies the
hypotheses of the Lemma, and the proof is completed.
\end{proof}

Case (5): $\left|  \mu^{\prime}\right|  =\left|  \lambda^{\prime}\right|
+1,\mu=\sigma-\varepsilon_{1}$ with $\sigma_{1}=\sigma_{2}$, is the only case
in the proof which depends on the $\mu^{\prime\prime}$ part of the formula (of
course the first steps of an induction proof are crucial!). By way of an
example, here is the matrix $B$ for $\mathcal{P}_{6}^{\left(  4\right)  }$%
\[
\left[
\begin{array}
[c]{rrrrrrrrr}%
1 & -2 & 0 & -2 & 8 & 0 & 2 & -18 & 18\\
0 & 1 & 0 & 0 & -10 & 0 & 0 & 42 & -72\\
0 & 0 & 1 & -3 & 3 & 0 & 3 & -9 & 9\\
0 & 0 & 0 & 1 & -3 & 0 & -2 & 13 & -18\\
0 & 0 & 0 & 0 & 1 & 0 & 0 & -9 & 24\\
0 & 0 & 0 & 0 & 0 & 1 & -2 & 2 & -2\\
0 & 0 & 0 & 0 & 0 & 0 & 1 & -4 & 9\\
0 & 0 & 0 & 0 & 0 & 0 & 0 & 1 & -6\\
0 & 0 & 0 & 0 & 0 & 0 & 0 & 0 & 1
\end{array}
\right]
\]
with the rows and columns labeled by $\left[
2211,3111,222,321,411,33,42,51,6\right]  $ (suppressing the commas in the partitions).

\section{Invariant harmonic polynomials}

In the sequel the symmetric polynomials have the argument $x^{2}=\left(
x_{1}^{2},\ldots,x_{N}^{2}\right)  $ (and $y^{2}$ is similarly defined). Thus
$e_{1}\left(  x^{2}\right)  =\left|  \left|  x\right|  \right|  ^{2}$, and the
former will often be used in equations involving symmetric functions. We
recall the basic facts about the Poisson kernel, that is, the reproducing
kernel, for harmonic polynomials $\mathbb{H}_{n}$, (see Dunkl and Xu \cite[Ch.
5]{DX}). The intertwining operator $V$ is a linear isomorphism on polynomials
such that $V\mathbb{P}_{n}^{\left(  N\right)  }=\mathbb{P}_{n}^{\left(
N\right)  }$ for each $n$, $\mathcal{D}_{i}V=V\frac{\partial}{\partial x_{i}}$
for $1\leq i\leq N$, and $V1=1$. Here the reflection group is a direct product
and so the intertwining map is the $N$-fold tensor product of the
one-dimensional transform (see \cite[Theorem 5.1]{D2}), defined by (for
$n\geq0$)
\[
Vx_{1}^{2n}=\frac{\left(  \frac{1}{2}\right)  _{n}}{\left(  \kappa+\frac{1}%
{2}\right)  _{n}}x_{1}^{2n},Vx_{1}^{2n+1}=\frac{\left(  \frac{1}{2}\right)
_{n+1}}{\left(  \kappa+\frac{1}{2}\right)  _{n+1}}x_{1}^{2n+1}.
\]
Then let $K_{n}\left(  x,y\right)  =\frac{1}{n!}V^{\left(  x\right)  }\left(
\left\langle x,y\right\rangle ^{n}\right)  $, for $x,y\in\mathbb{R}^{N}%
,n\geq0$ ($V^{\left(  x\right)  }$ acts on $x$). The key properties of $K_{n}$
needed here are $K_{n}\left(  xw,yw\right)  =K_{n}\left(  x,y\right)  $ for
$w\in B_{N}$ and $\mathcal{D}_{i}^{\left(  x\right)  }K_{n}\left(  x,y\right)
=y_{i}K_{n-1}\left(  x,y\right)  $ for $1\leq i\leq N$.

\begin{definition}
Let $d\omega$ denote the normalized rotation-invariant surface measure on the
unit sphere $S=\left\{  x\in\mathbb{R}^{N}:\left|  \left|  x\right|  \right|
=1\right\}  $, and for polynomials $f,g$ the inner product is
\[
\left\langle f,g\right\rangle _{S}=c_{\kappa}\int_{S}f\left(  x\right)
g\left(  x\right)  \prod_{i=1}^{N}\left|  x_{i}\right|  ^{2\kappa}%
\,d\omega\left(  x\right)  ,
\]
where $c_{\kappa}=\dfrac{\Gamma\left(  \frac{N}{2}+N\kappa\right)
\Gamma\left(  \frac{1}{2}\right)  ^{N}}{\Gamma\left(  \frac{N}{2}\right)
\Gamma\left(  \kappa+\frac{1}{2}\right)  ^{N}}$ so that $\left\langle
1,1\right\rangle _{S}=1$.
\end{definition}

The reproducing kernel for $\mathbb{H}_{n},n\geq1$ is
\begin{align*}
P_{n}\left(  x,y\right)   &  =\left(  \frac{N}{2}+N\kappa+n-1\right) \\
&  \times\sum_{j\leq n/2}\left(  -1\right)  ^{j}\frac{\left(  \frac{N}%
{2}+N\kappa\right)  _{n-1-j}}{j!}2^{n-2j}\left|  \left|  x\right|  \right|
^{2j}\left|  \left|  y\right|  \right|  ^{2j}K_{n-2j}\left(  x,y\right)  ,
\end{align*}
and satisfies $\left\langle P_{n}\left(  \cdot,y\right)  ,f\right\rangle
_{S}=f\left(  y\right)  $ for each $f\in\mathbb{H}_{n}$. The general theory
\cite{D1} for our differential-difference operators shows that $f\in
\mathbb{H}_{n},\allowbreak g\in\mathbb{H}_{m},\allowbreak m\neq n$ implies
$\left\langle f,g\right\rangle _{S}=0$. Let $K_{2n}^{0}$ denote the kernel
$K_{2n}\left(  x,y\right)  $ symmetrized with respect to the group
$\mathbb{Z}_{2}^{N}$ (for fixed $x\in\mathbb{R}^{N}$ sum over the $2^{N}$
points $\left(  \pm x_{1},\ldots,\pm x_{N}\right)  $ and divide by $2^{N}$);
it is the reproducing kernel for $\mathbb{H}_{2n}^{0}$. Let $K_{2n}^{B}\left(
x,y\right)  ,P_{2n}^{B}\left(  x,y\right)  $ denote the symmetrizations of
$K_{2n}\left(  x,y\right)  ,P_{2n}\left(  x,y\right)  $ respectively, with
respect to the group $B_{N}$. Thus $K_{2n}^{B}\left(  x,y\right)  =\frac
{1}{N!}\sum_{w\in S_{N}}K_{2n}^{0}\left(  xw,y\right)  $.

\begin{proposition}
For $n\geq1,\,x,y\in\mathbb{R}^{N}$%
\begin{align*}
K_{2n}^{0}\left(  x,y\right)   &  =2^{-2n}\sum_{\alpha\in\mathbb{N}_{0}%
^{N},\left|  \alpha\right|  =n}\frac{1}{\alpha!\,\left(  \kappa+\frac{1}%
{2}\right)  _{\alpha}}x^{2\alpha}y^{2\alpha},\\
K_{2n}^{B}\left(  x,y\right)   &  =2^{-2n}\sum_{\lambda\in\mathcal{P}%
_{n}^{\left(  N\right)  }}\frac{1}{\lambda!\,\left(  \kappa+\frac{1}%
{2}\right)  _{\lambda}}\frac{m_{\lambda}\left(  x^{2}\right)  m_{\lambda
}\left(  y^{2}\right)  }{m_{\lambda}\left(  1^{N}\right)  }.
\end{align*}
\end{proposition}

\begin{proof}
By the multinomial theorem, $\frac{\left\langle x,y\right\rangle ^{2n}%
}{\left(  2n\right)  !}=$\allowbreak$\frac{1}{\left(  2n\right)  !}%
\sum\left\{  \binom{2n}{\beta}x^{\beta}y^{\beta}:\beta\in\mathbb{N}_{0}%
^{N},\allowbreak\left|  \beta\right|  =2n\right\}  $. Symmetrizing over
$\mathbb{Z}_{2}^{N}$ removes monomials with odd exponents, and thus
$K_{2n}^{0}\left(  x,y\right)  $ is the result of applying $V$ to
$\sum\left\{  \frac{1}{\left(  2\alpha\right)  !}x^{2\alpha}y^{2\alpha}%
:\alpha\in\mathbb{N}_{0}^{N},\left|  \alpha\right|  =n\right\}  $. Further
$\left(  2\alpha\right)  !=2^{2\left|  \alpha\right|  }\alpha!\left(  \frac
{1}{2}\right)  _{\alpha}$, and $Vx^{2\alpha}=\left(  \left(  \frac{1}%
{2}\right)  _{\alpha}/\left(  \kappa+\frac{1}{2}\right)  _{\alpha}\right)
x^{2\alpha}$. For a fixed $\lambda\in\mathcal{P}_{n}^{\left(  N\right)  }$ the
sum \newline $\frac{1}{N!}\sum\left\{  x^{2w\alpha}y^{2\alpha}:w\in
S_{N},\alpha^{+}=\lambda\right\}  $ is a multiple of $m_{\lambda}\left(
x^{2}\right)  m_{\lambda}\left(  y^{2}\right)  $, and evaluated at $x=1^{N}=y$
the sum equals $m_{\lambda}\left(  1^{N}\right)  $ (that is, $\#\left\{
\alpha\in\mathbb{N}_{0}^{N}:\alpha^{+}=\lambda\right\}  $). The terms
$\alpha!$ and $\left(  \kappa+\frac{1}{2}\right)  _{\alpha}$ are invariant
under $S_{N}$.
\end{proof}

\begin{corollary}
For $n\geq1,\,x,y\in\mathbb{R}^{N}$
\begin{align*}
P_{2n}^{B}\left(  x,y\right)   &  =\left(  \frac{N}{2}+N\kappa+2n-1\right)
\sum_{j=0}^{n}\left(  -1\right)  ^{j}\frac{\left(  \frac{N}{2}+N\kappa\right)
_{2n-1-j}}{j!}e_{1}\left(  x^{2}\right)  ^{j}e_{1}\left(  y^{2}\right)  ^{j}\\
&  \times\sum_{\lambda\in\mathcal{P}_{n-j}^{\left(  N\right)  }}\frac
{1}{\lambda!\,\left(  \kappa+\frac{1}{2}\right)  _{\lambda}}\frac{m_{\lambda
}\left(  x^{2}\right)  m_{\lambda}\left(  y^{2}\right)  }{m_{\lambda}\left(
1^{N}\right)  }.
\end{align*}
\end{corollary}

Let $\mathbb{P}_{2n}^{B}$ denote the space of $B_{N}$-invariant elements of
$\mathbb{P}_{2n}^{\left(  N\right)  }$, then $\dim\mathbb{P}_{2n}%
^{B}=\#\mathcal{P}_{n}^{\left(  N\right)  }$. The space $\mathbb{H}_{2n}^{B}$
of $B_{N}$-invariant harmonic polynomials is the kernel of $\Delta_{\kappa}$;
since $\Delta_{\kappa}$ commutes with the action of $B_{N}$ and maps
$\mathbb{P}_{2n}^{\left(  N\right)  }$ onto $\mathbb{P}_{2n-2}^{\left(
N\right)  }$ we see that $\dim\mathbb{H}_{2n}^{B}=\#\mathcal{P}_{n}^{\left(
N\right)  }-\#\mathcal{P}_{n-1}^{\left(  N\right)  }$. The map $\lambda
\mapsto\lambda+\varepsilon_{1}$ is a one-to-one correspondence of
$\mathcal{P}_{n-1}^{\left(  N\right)  }$ onto a subset of $\mathcal{P}%
_{n}^{\left(  N\right)  }$ whose complement is
\[
\widetilde{\mathcal{P}}_{n}^{\left(  N\right)  }=\left\{  \lambda
\in\mathcal{P}_{n}^{\left(  N\right)  }:\lambda_{1}=\lambda_{2}\right\}  .
\]
Thus $\dim\mathbb{H}_{2n}^{B}=\#\widetilde{\mathcal{P}}_{n}^{\left(  N\right)
}$, and we will construct a basis for $\mathbb{H}_{2n}^{B}$ whose elements are
labeled in a natural way by $\widetilde{\mathcal{P}}_{n}^{\left(  N\right)  }%
$. There is a generating function for the dimensions: for $\lambda
\in\widetilde{\mathcal{P}}_{n}^{\left(  N\right)  }$ the conjugate
$\lambda^{T}$ (the partition corresponding to the transpose of the Ferrers
diagram) of $\lambda$ is a partition with $2\leq\left(  \lambda^{T}\right)
_{i}\leq N$ for all $i$ (no parts equal to 1 or exceeding $N$); thus
\[
\sum_{n=0}^{\infty}\left(  \#\widetilde{\mathcal{P}}_{n}^{\left(  N\right)
}\right)  q^{n}=\prod_{j=2}^{N}\left(  1-q^{j}\right)  ^{-1}.
\]
This expression yields an estimate for $\#\widetilde{\mathcal{P}}_{n}^{\left(
N\right)  }$. Indeed, let $M=\operatorname{lcm}\left(  2,3,\ldots,N\right)  $,
$\prod_{j=2}^{N}\left(  1-q^{j}\right)  ^{-1}=p\left(  q\right)  \left(
1-q^{M}\right)  ^{-\left(  N-1\right)  }$ for some polynomial $p\left(
q\right)  $; thus $\#\widetilde{\mathcal{P}}_{n}^{\left(  N\right)  }=O\left(
\left(  \frac{n}{M}\right)  ^{N-2}\right)  $ as $n\rightarrow\infty$. In
$K_{2n}^{B}\left(  x,y\right)  $ expand each $m_{\lambda}\left(  y^{2}\right)
$ in the $\widetilde{m}_{\mu}$-basis
\[
K_{2n}^{B}\left(  x,y\right)  =2^{-2n}\sum_{\mu\preceq\lambda}\frac{B\left(
\mu,\lambda\right)  }{\lambda!\,\left(  \kappa+\frac{1}{2}\right)  _{\lambda}%
}\frac{m_{\lambda}\left(  x^{2}\right)  }{m_{\lambda}\left(  1^{N}\right)
}\widetilde{m}_{\mu}\left(  y^{2}\right)  ;
\]
this leads to the following.

\begin{definition}
For $\mu\in\mathcal{P}_{n}^{\left(  N\right)  }$ the $B_{N}$-invariant
polynomial $h_{\mu}$ is given by
\[
h_{\mu}\left(  x\right)  =2^{-n}\sum_{\lambda\in\mathcal{P}_{n}^{\left(
N\right)  },\mu\preceq\lambda}\frac{B\left(  \mu,\lambda\right)  }%
{\lambda!\,\left(  \kappa+\frac{1}{2}\right)  _{\lambda}}\frac{m_{\lambda
}\left(  x^{2}\right)  }{m_{\lambda}\left(  1^{N}\right)  }.
\]
\end{definition}

Recall that $\lambda\in\mathcal{P}_{n}$ and $B\left(  \lambda,\mu\right)
\neq0$ implies $l\left(  \lambda\right)  \leq l\left(  \mu\right)  $ and
$\lambda\in\mathcal{P}_{n}^{\left(  N\right)  }$.

\begin{theorem}
\label{laph}For $\mu\in\mathcal{P}_{n}^{\left(  N\right)  }$, if $\mu_{1}%
>\mu_{2}$ then $\Delta_{\kappa}h_{\mu}=2h_{\mu-\varepsilon_{1}}$ and if
$\mu_{1}=\mu_{2}$ then $\Delta_{\kappa}h_{\mu}=0$.
\end{theorem}

\begin{proof}
From the basic properties of $K_{2n}$ it follows that $\Delta_{\kappa
}^{\left(  x\right)  }K_{2n}\left(  x,y\right)  =\left|  \left|  y\right|
\right|  ^{2}K_{2n-2}\left(  x,y\right)  $. Symmetrize this equation with
respect to $B_{N}$ (and note that $\Delta_{\kappa}$ commutes with the group
action) to obtain
\begin{align*}
\Delta_{\kappa}^{\left(  x\right)  }K_{2n}^{B}\left(  x,y\right)   &
=2^{-n}\sum_{\mu\in\mathcal{P}_{n}^{\left(  N\right)  }}\Delta_{\kappa
}^{\left(  x\right)  }h_{\mu}\left(  x\right)  \widetilde{m}_{\mu}\left(
y^{2}\right) \\
&  =e_{1}\left(  y^{2}\right)  K_{2n-2}^{B}\left(  x,y\right) \\
&  =2^{1-n}\sum_{\sigma\in\mathcal{P}_{n-1}^{\left(  N\right)  }}h_{\sigma
}\left(  x\right)  \widetilde{m}_{\sigma+\varepsilon_{1}}\left(  y^{2}\right)
.
\end{align*}
By definition $e_{1}\left(  y^{2}\right)  \widetilde{m}_{\sigma}\left(
y^{2}\right)  =\widetilde{m}_{\sigma+\varepsilon_{1}}\left(  y^{2}\right)  $.
Considering the equations as expansions in $\left\{  \widetilde{m}_{\mu
}\left(  y^{2}\right)  :\mu\in\mathcal{P}_{n}^{\left(  N\right)  }\right\}  $
shows $\Delta_{\kappa}h_{\sigma+\varepsilon_{1}}\left(  x\right)  =2h_{\sigma
}\left(  x\right)  $ for $\sigma\in\mathcal{P}_{n-1}^{\left(  N\right)  }$ and
$\Delta_{\kappa}h_{\mu}\left(  x\right)  =0$ if $\mu_{1}=\mu_{2}$.
\end{proof}

\begin{corollary}
The set $\left\{  h_{\mu}:\mu\in\widetilde{\mathcal{P}}_{n}^{\left(  N\right)
}\right\}  $ is a basis for $\mathbb{H}_{2n}^{B}$.
\end{corollary}

\begin{proof}
For $\lambda,\mu\in\widetilde{\mathcal{P}}_{n}^{\left(  N\right)  }$ the
coefficient of $m_{\lambda}$ in $h_{\mu}$ is $\delta_{\lambda\mu}\left(
2^{n}\lambda!\,\left(  \kappa+\frac{1}{2}\right)  _{\lambda}m_{\lambda}\left(
1^{N}\right)  \right)  ^{-1}$ and thus $\left\{  h_{\mu}:\mu\in\widetilde
{\mathcal{P}}_{n}^{\left(  N\right)  }\right\}  $ is linearly independent.
Also $\dim\mathbb{H}_{2n}^{B}=\#\widetilde{\mathcal{P}}_{n}^{\left(  N\right)
}$.
\end{proof}

Besides the inner product on polynomials defined by integration over the
sphere $S$ there is also the important pairing defined in an algebraic manner
using the operators $\mathcal{D}_{i}$, namely
\[
\left\langle f,g\right\rangle _{h}=f\left(  \mathcal{D}_{1},\ldots
,\mathcal{D}_{N}\right)  g\left(  x\right)  |_{x=0}.
\]
Since $\mathcal{D}_{i}^{2}x_{i}^{2n}=2n\left(  2n-1+2\kappa\right)
x_{i}^{2n-2}$ and $\mathcal{D}_{j}^{2}x_{i}^{2n}=0$ for $j\neq i,$ we have
that $\left\langle x_{i}^{2n},x_{i}^{2n}\right\rangle _{h}=2^{2n}n!\left(
\kappa+\frac{1}{2}\right)  _{n}$ and $\left\langle x^{2\alpha},x^{2\beta
}\right\rangle _{h}=\delta_{\alpha\beta}2^{2\left|  \alpha\right|  }%
\alpha!\left(  \kappa+\frac{1}{2}\right)  _{\alpha}$ for $\alpha,\beta
\in\mathbb{N}_{0}^{N}$, so that monomials are mutually orthogonal. It follows
that $\left\langle m_{\lambda}\left(  x^{2}\right)  ,m_{\mu}\left(
x^{2}\right)  \right\rangle _{h}=\delta_{\lambda\mu}2^{2n}\lambda!\left(
\kappa+\frac{1}{2}\right)  _{\lambda}m_{\lambda}\left(  1^{N}\right)  $ for
$\lambda,\mu\in\mathcal{P}_{n}^{\left(  N\right)  }$. It was shown in (
\cite[Theorem 3.8]{D2}, see also \cite[Theorem 5.2.4]{DX}) that for
$f,g\in\mathbb{H}_{n}$
\[
\left\langle f,g\right\rangle _{h}=2^{n}\left(  \frac{N}{2}+N\kappa\right)
_{n}\left\langle f,g\right\rangle _{S}.
\]
This allows the direct calculation of $\left\langle h_{\lambda},h_{\mu
}\right\rangle _{S}$ for $\lambda,\mu\in\widetilde{\mathcal{P}}_{n}^{\left(
N\right)  }$.

\begin{proposition}
\label{inpro1}For $\lambda,\mu\in\widetilde{\mathcal{P}}_{n}^{\left(
N\right)  }$
\begin{align*}
\left\langle h_{\lambda},h_{\mu}\right\rangle _{h}  &  =\sum_{\lambda
\preceq\sigma,\mu\preceq\sigma}\frac{B\left(  \lambda,\sigma\right)  B\left(
\mu,\sigma\right)  }{\sigma!\left(  \kappa+\frac{1}{2}\right)  _{\sigma
}m_{\sigma}\left(  1^{N}\right)  },\\
\left\langle h_{\lambda},h_{\mu}\right\rangle _{S}  &  =2^{-2n}\left(  \left(
\frac{N}{2}+N\kappa\right)  _{2n}\right)  ^{-1}\left\langle h_{\lambda}%
,h_{\mu}\right\rangle _{h}.
\end{align*}
\end{proposition}

Computations with small degrees ($n\leq8$) show that there are no product
(linear factors in $\kappa$) formulae for the inner products and $\left\{
h_{\lambda}\right\}  $ is not an orthogonal basis. It is however possible to
find a nice formula for the determinant of the Gram matrix $\left(
\left\langle h_{\lambda},h_{\mu}\right\rangle _{h}\right)  _{\lambda,\mu
\in\widetilde{\mathcal{P}}_{n}^{\left(  N\right)  }}$. Before presenting this
result we consider the bi-orthogonal set for $\left\{  h_{\mu}:\mu
\in\widetilde{\mathcal{P}}_{n}^{\left(  N\right)  }\right\}  $ constructed by
extracting a certain multiple of the coefficient of $m_{\lambda}\left(
y^{2}\right)  $ for $\lambda\in\widetilde{\mathcal{P}}_{n}^{\left(  N\right)
}$ in the Poisson kernel $P_{2n}^{B}\left(  x,y\right)  $. Indeed $P_{2n}^{B}$
is a multiple of
\begin{align*}
&  \sum_{j=0}^{n}\sum_{\sigma\in\mathcal{P}_{n-j}^{\left(  N\right)  }}\left(
-1\right)  ^{j}e_{1}\left(  x^{2}\right)  ^{j}e_{1}\left(  y^{2}\right)
^{j}\frac{\left(  \frac{N}{2}+N\kappa\right)  _{2n-1-j}}{j!\,\sigma!\,\left(
\kappa+\frac{1}{2}\right)  _{\sigma}}\frac{m_{\sigma}\left(  x^{2}\right)
m_{\sigma}\left(  y^{2}\right)  }{m_{\sigma}\left(  1^{N}\right)  }\\
&  =\sum_{j=0}^{n}\sum_{\lambda\in\mathcal{P}_{n}^{\left(  N\right)  }}%
\sum_{\sigma\in\mathcal{P}_{n-j}^{\left(  N\right)  }}\left(  -1\right)
^{j}e_{1}\left(  x^{2}\right)  ^{j}%
\genfrac{\langle}{\rangle}{0pt}{}{\lambda}{\sigma}%
\frac{\left(  \frac{N}{2}+N\kappa\right)  _{2n-1-j}}{j!\,\sigma!\,\left(
\kappa+\frac{1}{2}\right)  _{\sigma}}\frac{m_{\sigma}\left(  x^{2}\right)
m_{\lambda}\left(  y^{2}\right)  }{m_{\sigma}\left(  1^{N}\right)  }.
\end{align*}
Then let
\[
g_{\lambda}\left(  x\right)  =\sum_{j=0}^{n}\frac{e_{1}\left(  x^{2}\right)
^{j}}{j!\,\left(  -\frac{N}{2}-N\kappa-2n+2\right)  _{j}}\sum_{\sigma
\in\mathcal{P}_{n-j}^{\left(  N\right)  },\sigma\subset\lambda}%
\genfrac{\langle}{\rangle}{0pt}{}{\lambda}{\sigma}%
\frac{1}{\sigma!\,\left(  \kappa+\frac{1}{2}\right)  _{\sigma}}\frac
{m_{\sigma}\left(  x^{2}\right)  }{m_{\sigma}\left(  1^{N}\right)  },
\]
so that the coefficient of $m_{\lambda}\left(  y^{2}\right)  $ in $P_{2n}%
^{B}\left(  x,y\right)  $ is $\left(  \frac{N}{2}+N\kappa\right)
_{2n}g_{\lambda}\left(  x\right)  $. Observe that $g_{\lambda}\left(
x\right)  =\dfrac{1}{\lambda!\,\left(  \kappa+\frac{1}{2}\right)  _{\lambda}%
}\dfrac{m_{\lambda}\left(  x^{2}\right)  }{m_{\lambda}\left(  1^{N}\right)
}+\allowbreak\left|  \left|  x\right|  \right|  ^{2}g_{\lambda}^{\prime
}\left(  x\right)  $ where $g_{\lambda}^{\prime}\left(  x\right)  $ is a
polynomial of degree $2n-2$.

\begin{proposition}
For $\lambda,\mu\in\widetilde{\mathcal{P}}_{n}^{\left(  N\right)  }$ the inner
product
\[
\left\langle g_{\lambda},h_{\mu}\right\rangle _{h}=\delta_{\lambda\mu}%
\frac{2^{n}}{\lambda!\,\left(  \kappa+\frac{1}{2}\right)  _{\lambda}%
m_{\lambda}\left(  1^{N}\right)  }.
\]
\end{proposition}

\begin{proof}
Indeed $\left\langle g_{\lambda},h_{\mu}\right\rangle _{h}=\frac{1}%
{\lambda!\,\left(  \kappa+\frac{1}{2}\right)  _{\lambda}m_{\lambda}\left(
1^{N}\right)  }\left\langle m_{\lambda}\left(  x^{2}\right)  ,h_{\mu}\left(
x\right)  \right\rangle _{h}+\allowbreak g_{\lambda}^{\prime}\left(
\mathcal{D}\right)  \Delta_{\kappa}h_{\mu}\left(  x\right)  =\allowbreak
\delta_{\lambda\mu}\frac{2^{n}}{\lambda!\,\left(  \kappa+\frac{1}{2}\right)
_{\lambda}m_{\lambda}\left(  1^{N}\right)  }$, by the properties of the
pairing $\left\langle \cdot,\cdot\right\rangle _{h}$.
\end{proof}

Essentially $g_{\lambda}$ is (a scalar multiple of) the projection of
$m_{\lambda}\left(  x^{2}\right)  $ onto $\mathbb{H}_{2n}$. The method of
orthogonal projections to construct bases of harmonic polynomials was studied
in more detail by Y. Xu \cite{X1}.

We turn to the evaluation of the determinant of the Gram matrix of
$\{h_{\lambda}:\lambda\in\widetilde{\mathcal{P}}_{n}^{\left(  N\right)  }\}$
using the $\left\langle \cdot,\cdot\right\rangle _{h}$ inner product; the
value for $\left\langle \cdot,\cdot\right\rangle _{S}$ is just a product of a
power ($\#\widetilde{\mathcal{P}}_{n}^{\left(  N\right)  }$) of the
proportionality factor with the previous one. It is perhaps surprising that
the calculation can be carried out by finding the determinant for the entire
set $\left\{  h_{\lambda}:\lambda\in\mathcal{P}_{n}^{\left(  N\right)
}\right\}  $, which can easily be done. By the orthogonality of $\left\{
m_{\mu}\right\}  $ we see that the Gram matrix $G\left(  \lambda,\mu\right)
=\left\langle h_{\lambda},h_{\mu}\right\rangle _{h}$ for $\lambda,\mu
\in\mathcal{P}_{n}^{\left(  N\right)  }$ satisfies
\begin{align*}
G &  =BCB^{T},\\
\det G &  =\left(  \det B\right)  ^{2}\det C=\det C,
\end{align*}
where the diagonal matrix $C\left(  \lambda,\mu\right)  =\delta_{\lambda\mu
}\left(  \lambda!\,\left(  \kappa+\frac{1}{2}\right)  _{\lambda}\,m_{\lambda
}\left(  1^{N}\right)  \right)  ^{-1}$; just as in Proposition \ref{inpro1}
$\left\langle h_{\lambda},h_{\mu}\right\rangle _{h}=\sum\limits_{\lambda
\preceq\sigma,\mu\preceq\sigma}\frac{B\left(  \lambda,\sigma\right)  B\left(
\mu,\sigma\right)  }{\sigma!\left(  \kappa+\frac{1}{2}\right)  _{\sigma
}m_{\sigma}\left(  1^{N}\right)  }$ for any $\lambda,\mu\in\mathcal{P}%
_{n}^{\left(  N\right)  }$. We already showed $B$ is triangular with 1's on
the main diagonal so that $\det B=1$. Let
\[
D_{n}=\det C=\prod_{\lambda\in\mathcal{P}_{n}^{\left(  N\right)  }}\left(
\lambda!\,\left(  \kappa+\frac{1}{2}\right)  _{\lambda}\,m_{\lambda}\left(
1^{N}\right)  \right)  ^{-1},
\]
and let $\widetilde{D}_{n}=\det\widetilde{G}$, where $\widetilde{G}$ is the
principal submatrix of $G$ for the labels $\lambda,\mu\in\widetilde
{\mathcal{P}}_{n}^{\left(  N\right)  }$. We will show that
\[
\widetilde{D}_{n}=\frac{D_{n}}{D_{n-1}}\prod_{\lambda\in\mathcal{P}%
_{n-1}^{\left(  N\right)  }}\left(  \left(  \lambda_{1}-\lambda_{2}+1\right)
\left(  N\kappa+\frac{N}{2}+2\left(  n-1)-(\lambda_{1}-\lambda_{2}\right)
\right)  \right)  .
\]
A simplified form of this will be given later. The idea of the proof is to use
the orthogonal decomposition $\mathbb{P}_{n}^{\left(  N\right)  }%
=\bigoplus_{j\leq n/2}\left|  \left|  x\right|  \right|  ^{2j}\mathbb{H}%
_{n-2j}$ to produce a transformation of the Gram matrix into block form. The
underlying relation is the product formula for $\Delta_{\kappa}$: for
$m,j\geq0$ and any $f\in\mathbb{P}_{m}^{\left(  N\right)  }$%
\[
\Delta_{\kappa}\left|  \left|  x\right|  \right|  ^{2j}f\left(  x\right)
=4j\left(  m+j-1+\frac{N}{2}+N\kappa\right)  \left|  \left|  x\right|
\right|  ^{2j-2}f\left(  x\right)  +\left|  \left|  x\right|  \right|
^{2j}\Delta_{\kappa}f\left(  x\right)  .
\]
Specializing to harmonic polynomials $f\in\mathbb{H}_{m}$ and iterating this
formula shows that
\begin{equation}
\Delta_{\kappa}^{s}\left|  \left|  x\right|  \right|  ^{2j}f\left(  x\right)
=2^{2s}\left(  -j\right)  _{s}\left(  -m-j+1-\frac{N}{2}-N\kappa\right)
_{s}\left|  \left|  x\right|  \right|  ^{2j-2s}f\left(  x\right)
\label{dpower}%
\end{equation}
for $s\geq0$; note $\Delta_{\kappa}^{s}\left|  \left|  x\right|  \right|
^{2j}f\left(  x\right)  =0$ for $s>j$ (see \cite[Theorem 3.6]{D2}). Suppose
$f\in\mathbb{H}_{n-2i},g\in\mathbb{H}_{n-2j}$ ($i,j\leq\frac{n}{2}$) then
\[
\left\langle \left|  \left|  x\right|  \right|  ^{2i}f\left(  x\right)
,\left|  \left|  x\right|  \right|  ^{2j}g\left(  x\right)  \right\rangle
_{h}=\delta_{ij}2^{2i}i!\left(  n-2i+\frac{N}{2}+N\kappa\right)
_{i}\left\langle f,g\right\rangle _{h},
\]
(by the symmetry of the inner product we can assume $i>j$ and use equation
\ref{dpower}).

\begin{lemma}
\label{dpow1}Suppose $f\in\mathbb{P}_{n}^{\left(  N\right)  }$ and $f\left(
x\right)  =\sum_{j\leq n/2}\left|  \left|  x\right|  \right|  ^{2j}%
f_{n-2j}\left(  x\right)  $ with $f_{n-2j}\in\mathbb{H}_{n-2j}$ for each
$j\leq\frac{n}{2}$ then $\Delta_{\kappa}^{i}f=2^{2i}i!\left(  n-2i+\frac{N}%
{2}+N\kappa\right)  _{i}\,f_{n-2i}$ if and only if $f_{n-2j}=0$ for each
$j>i$, (that is, $\Delta_{\kappa}^{i+1}f=0$).
\end{lemma}

\begin{proof}
By equation \ref{dpower} and for each $i\leq\frac{n}{2}$
\[
\Delta_{\kappa}^{i}f\left(  x\right)  =\sum_{i\leq j\leq n/2}2^{2i}\left(
-j\right)  _{i}\left(  -n+j+1-\frac{N}{2}-N\kappa\right)  _{i}\left|  \left|
x\right|  \right|  ^{2j-2i}f_{n-2j}\left(  x\right)  .
\]
This is an orthogonal expansion, hence $\Delta_{\kappa}^{i}f\left(  x\right)
$ equals the term in the sum with $j=i$ if and only if $f_{n-2j}=0$ for all
$j>i$.
\end{proof}

The Lemma allows us to find the lowest-degree part of the expansion of
$h_{\mu}$ for $\mu\in\mathcal{P}_{n}^{\left(  N\right)  }$, since
$\Delta_{\kappa}^{i}h_{\mu}=2^{i}h_{\mu-i\varepsilon_{1}}$ if $i\leq\mu
_{1}-\mu_{2}$ and $\Delta_{\kappa}^{i}h_{\mu}=0$ if $i>\mu_{1}-\mu_{2}$,
by\ Theorem \ref{laph}. We will use an elementary property of any inner
product space $E$: suppose $\left\{  f_{i}:1\leq i\leq m\right\}  $ is a
linearly independent set in $E$, for some $k<m$ let $\rho$ denote the
orthogonal projection of $E$ onto $\mathrm{span}\left\{  f_{i}:1\leq i\leq
k\right\}  $, then
\[
\det\left(  \left\langle f_{i},f_{j}\right\rangle \right)  _{i,j=1}^{m}%
=\det\left(  \left\langle f_{i},f_{j}\right\rangle \right)  _{i,j=1}^{k}%
\det\left(  \left\langle \left(  1-\rho\right)  f_{i},\left(  1-\rho\right)
f_{j}\right\rangle \right)  _{i,j=k+1}^{m}.
\]
The proof is easy: for each $i>k$ there are coefficients $a_{ji}$ for $1\leq
j\leq m$ so that $\rho f_{i}=\sum_{j=1}^{k}a_{ji}f_{j}$. Let $A$ be the
$m\times m$ matrix with entries $A_{ji}=\delta_{ji}$ except $A_{ji}=-a_{ji}$
for $1\leq j\leq k<i\leq m$. Let $G$ be the Gram matrix of $\left\{
f_{i}\right\}  _{i=1}^{m}$ and let $G^{\prime}=A^{T}GA$ so that
\[
G_{ij}^{\prime}=G_{ji}^{\prime}=\left\{
\begin{tabular}
[c]{l}%
$\left\langle f_{i},f_{j}\right\rangle $ for $1\leq i\leq j\leq k$\\
$\left\langle f_{i},\left(  1-\rho\right)  f_{j}\right\rangle =0$ for $1\leq
i\leq k<j\leq m$\\
$\left\langle \left(  1-\rho\right)  f_{i},\left(  1-\rho\right)
f_{j}\right\rangle $ for $k+1\leq i\leq j\leq m;$%
\end{tabular}
\right.
\]
thus $\det G^{\prime}=\left(  \det A\right)  ^{2}\det G=\det G$ and
$G^{\prime}$ has a $2\times2$ block structure with $0$ in the off-diagonal
blocks. We apply this to the projection of $\mathbb{P}_{2n}^{B}$ onto
$\mathbb{H}_{2n}^{B}$, denoted by $\rho_{n}$. We write $e_{1}$ for $\left|
\left|  x\right|  \right|  ^{2}$ as before. For $\mu\in\mathcal{P}%
_{n}^{\left(  N\right)  }$ let $h_{\mu}=\sum_{j=0}^{\mu_{1}-\mu_{2}}e_{1}%
^{j}h_{\mu,j}\left(  x\right)  $ with each $h_{\mu,j}\in\mathbb{H}_{2n-2j}%
^{B}$ (and by the above discussion, $h_{\mu,i}$ is a multiple of
$h_{\mu-i\varepsilon_{1}}$ for $i=\mu_{1}-\mu_{2}$), then $\left(  1-\rho
_{n}\right)  h_{\mu}=\sum_{j=1}^{\mu_{1}-\mu_{2}}e_{1}^{j}h_{\mu,j}$. Thus
$D_{n}=\widetilde{D}_{n}\det M$ where $M$ is the Gram matrix for $\left\{
\left(  1-\rho_{n}\right)  h_{\mu}:\mu\in\mathcal{P}_{n}^{\left(  N\right)
},\mu_{1}>\mu_{2}\right\}  $. The span of this set is $e_{1}\mathbb{P}%
_{2n-2}^{B}$, \ and $\mathbb{P}_{2n-2}^{B}$ can be decomposed just like
$\mathbb{P}_{2n}^{B}$. Indeed $\left\{  \left(  1-\rho_{n}\right)  h_{\mu}%
:\mu\in\mathcal{P}_{n}^{\left(  N\right)  },\allowbreak\mu_{1}=\mu
_{2}+1\right\}  $ is a basis for $e_{1}\mathbb{H}_{2n-2}^{B}$, since $\left(
1-\rho_{n}\right)  h_{\mu}=e_{1}h_{\mu,1}$, a nonzero multiple of $e_{1}%
h_{\mu-\varepsilon_{1}}$, for $\mu_{1}=\mu_{2}+1$. Repeat the previous step
with the projection $e_{1}\rho_{n-1}e_{1}^{-1}$ to express $\det M$ as a
product. Here is a formal inductive argument.

Let $M^{\left(  j\right)  }\left(  \lambda,\mu\right)  $ be the symmetric
matrix indexed by $\mathcal{P}_{n}^{\left(  N\right)  }$ with entries
\[
M^{\left(  j\right)  }\left(  \lambda,\mu\right)  =\left\{
\begin{tabular}
[c]{l}%
$\left\langle e_{1}^{i}h_{\lambda,i},e_{1}^{i}h_{\mu,i}\right\rangle _{h}$ for
$\lambda_{1}-\lambda_{2}=\mu_{1}-\mu_{2}=i<j$\\
$0$ for $\lambda_{1}-\lambda_{2}<\min\left(  \mu_{1}-\mu_{2},j\right)  $\\
$\left\langle \sum\limits_{i=j}^{\lambda_{1}-\lambda_{2}}e_{1}^{i}%
h_{\lambda,i},\sum\limits_{i=j}^{\mu_{1}-\mu_{2}}e_{1}^{i}h_{\mu
,i}\right\rangle _{h}$ for $\lambda_{1}-\lambda_{2},\mu_{1}-\mu_{2}\geq j.$%
\end{tabular}
\right.
\]
By the symmetry $M^{\left(  j\right)  }\left(  \lambda,\mu\right)  =0$ if
$\mu_{1}-\mu_{2}<\min\left(  \lambda_{1}-\lambda_{2},j\right)  $. The matrix
$M^{\left(  n\right)  }$ has a diagonal block decomposition (zero blocks off
the diagonal) with one block for each set of labels $\left\{  \lambda
\in\mathcal{P}_{n}^{\left(  N\right)  }:\lambda_{1}=\lambda_{2}+i\right\}  $,
equivalently $\widetilde{\mathcal{P}}_{n-i}^{\left(  N\right)  }$, for each
$i=0,1,\ldots,n-2,n$ (the set $\widetilde{\mathcal{P}}_{1}^{\left(  N\right)
} $ is empty). The matrix $M^{\left(  0\right)  }=G$, the Gram matrix of
$\left\{  h_{\lambda}:\lambda\in\mathcal{P}_{n}^{\left(  N\right)  }\right\}
$. We show that $\det G=\det M^{\left(  0\right)  }$ by proving $\det
M^{\left(  j\right)  }=\det M^{\left(  j+1\right)  }$ for each $j<n$. Only the
principal submatrix of $M^{\left(  j\right)  }$ labeled by $\lambda$ with
$\lambda_{1}-\lambda_{2}\geq j$ need be considered. This is the Gram matrix of
a certain basis for $e_{1}^{j}\mathbb{P}_{2n-2j}^{B}$; the projection
$e_{1}^{j}\rho_{n-j}e_{1}^{-j}$ maps this space onto $e_{1}^{j}\mathbb{H}%
_{2n-2j}^{B}=\mathrm{span}\left\{  e_{1}^{j}h_{\lambda,j}:\lambda_{1}%
=\lambda_{2}+j\right\}  $. Now observe that $\left(  1-e_{1}^{j}\rho
_{n-j}e_{1}^{-j}\right)  \sum\limits_{i=j}^{\lambda_{1}-\lambda_{2}}e_{1}%
^{i}h_{\lambda,i}=\sum\limits_{i=j+1}^{\lambda_{1}-\lambda_{2}}e_{1}%
^{i}h_{\lambda,i}$ for $\lambda$ with $\lambda_{1}-\lambda_{2}>j$. The
projection argument shows $\det M^{\left(  j\right)  }=\det M^{\left(
j+1\right)  }$.

The determinant of the principal submatrix of $M^{\left(  n\right)  }$ labeled
by \newline $\left\{  \lambda\in\mathcal{P}_{n}^{\left(  N\right)  }%
:\lambda_{1}=\lambda_{2}+i\right\}  $ is a multiple of $\widetilde{D}_{n-i}$.
Let
\[
c_{i}=\left(  2^{2i}i!\left(  2n-2i+\frac{N}{2}+N\kappa\right)  _{i}\right)
^{-1};
\]
by Lemma \ref{dpow1} $h_{\lambda,i}=c_{i}\Delta_{\kappa}^{i}h_{\lambda}%
=2^{i}c_{i}h_{\lambda-i\varepsilon_{1}}$ for $\lambda_{1}=\lambda_{2}+i$.
Further
\begin{align*}
\left\langle e_{1}^{i}h_{\lambda,i},e_{1}^{i}h_{\mu,i}\right\rangle _{h}  &
=c_{i}^{-1}\left\langle h_{\lambda,i},h_{\mu,i}\right\rangle _{h}\\
&  =2^{2i}c_{i}\,\left\langle h_{\lambda-i\varepsilon_{1}},h_{\mu
-i\varepsilon_{1}}\right\rangle _{h}%
\end{align*}
for $\lambda_{1}-\lambda_{2}=\mu_{1}-\mu_{2}=i$. The correspondence
$\sigma\mapsto\sigma+i\varepsilon_{1}$ is one-to-one from $\widetilde
{\mathcal{P}}_{n-i}^{\left(  N\right)  }$ to $\left\{  \lambda\in
\mathcal{P}_{n}^{\left(  N\right)  }:\lambda_{1}=\lambda_{2}+i\right\}  $.
Thus
\[
D_{n}=\det G=\prod_{i=0}^{n}\left(  i!\left(  2n-2i+\frac{N}{2}+N\kappa
\right)  _{i}\right)  ^{-d\left(  n-i,N\right)  }\widetilde{D}_{n-i},
\]
where $d\left(  j,N\right)  =\#\widetilde{\mathcal{P}}_{j}^{\left(  N\right)
}$, for $j\geq0$. In particular, $d\left(  1,N\right)  =0$ and $\widetilde
{D}_{1}=1$.

\begin{theorem}
The determinant $\widetilde{D}_{n}$ of the Gram matrix of $\left\{
h_{\lambda}:\lambda\in\widetilde{\mathcal{P}}_{n}^{\left(  N\right)
}\right\}  $ for the inner product $\left\langle \cdot,\cdot\right\rangle
_{h}$ satisfies
\[
\widetilde{D}_{n}=\frac{D_{n}}{D_{n-1}}\prod_{i=1}^{n}\left(  i\left(
2n-i-1+\frac{N}{2}+N\kappa\right)  \right)  ^{d\left(  n-i,N\right)  }.
\]
\end{theorem}

\begin{proof}
From the above result
\[
\frac{D_{n}}{D_{n-1}}=\widetilde{D}_{n}\frac{\prod_{j=0}^{n-1}\left(
j!\left(  2n-2-2j+\frac{N}{2}+N\kappa\right)  _{j}\right)  ^{d\left(
n-1-j,N\right)  }}{\prod_{i=1}^{n}\left(  i!\left(  2n-2i+\frac{N}{2}%
+N\kappa\right)  _{i}\right)  ^{d\left(  n-i,N\right)  }}.
\]
In the numerator replace $j$ by $i-1$; the ratio $\dfrac{\left(
2n-2i+\frac{N}{2}+N\kappa\right)  _{i-1}}{\left(  2n-2i+\frac{N}{2}%
+N\kappa\right)  _{i}}=$\newline $\dfrac{1}{2n-i-1+\frac{N}{2}+N\kappa}$, and
this completes the proof.
\end{proof}

In the ratio $\frac{D_{n}}{D_{n-1}}$ there is some cancellation by use of the
decomposition $\mathcal{P}_{n}^{\left(  N\right)  }=\widetilde{\mathcal{P}%
}_{n}^{\left(  N\right)  }\cup\left\{  \mu+\varepsilon_{1}:\mu\in
\mathcal{P}_{n-1}^{\left(  N\right)  }\right\}  $; indeed each $\mu
\in\mathcal{P}_{n-1}^{\left(  N\right)  }$ contributes
\[
\dfrac{\mu!\left(  \kappa+\frac{1}{2}\right)  _{\mu}m_{\mu}\left(
1^{N}\right)  }{\left(  \mu+\varepsilon_{1}\right)  !\left(  \kappa+\frac
{1}{2}\right)  _{\mu+\varepsilon_{1}}m_{\mu+\varepsilon_{1}}\left(
1^{N}\right)  }=\dfrac{1}{\left(  \mu_{1}+1\right)  \left(  \kappa+\frac{1}%
{2}+\mu_{1}\right)  \#\left\{  j:\mu_{j}=\mu_{1}\right\}  }.
\]
Note that $m_{\mu}\left(  1^{N}\right)  =N!/\prod_{s\geq0}\left(  \#\left\{
j:\mu_{j}=s\right\}  \right)  !$ so the change from $m_{\mu}\left(
1^{N}\right)  $ to $m_{\mu+\varepsilon_{1}}\left(  1^{N}\right)  $ is the
replacement of $s!$ by $1!\left(  s-1\right)  !$ where $s=\#\left\{  j:\mu
_{j}=\mu_{1}\right\}  $ (except when $\mu=\left(  0\right)  $ in which case
$m_{\mu}\left(  1^{N}\right)  /m_{\mu+\varepsilon_{1}}\left(  1^{N}\right)
=\frac{1}{N}$; this only affects the vacuous equation $\widetilde{D}_{1}=1$).
The other part of the expression for $\widetilde{D}_{n}$ is a product with
$\sum_{i=1}^{n}d\left(  n-i,N\right)  =\#\mathcal{P}_{n-1}^{\left(  N\right)
}$ terms (and if $\mu\in\mathcal{P}_{n-1}^{\left(  N\right)  }$ with $\mu
_{1}-\mu_{2}+1=i$ then $\mu\in\widetilde{\mathcal{P}}_{n-i}^{\left(  N\right)
}+\left(  i-1\right)  \varepsilon_{1},1\leq i\leq n$). This is the simplified
formula (for $n\geq2$):
\begin{align*}
\widetilde{D}_{n}  &  =\prod_{\lambda\in\widetilde{\mathcal{P}}_{n}^{\left(
N\right)  }}\left(  \lambda!\,\left(  \kappa+\frac{1}{2}\right)  _{\lambda
}m_{\lambda}\left(  1^{N}\right)  \right)  ^{-1}\\
&  \times\prod_{\mu\in\mathcal{P}_{n-1}^{\left(  N\right)  }}\frac{\left(
\mu_{1}-\mu_{2}+1\right)  \left(  2n-2-\mu_{1}+\mu_{2}+\frac{N}{2}%
+N\kappa\right)  }{\left(  \mu_{1}+1\right)  \left(  \kappa+\frac{1}{2}%
+\mu_{1}\right)  \#\left\{  j:\mu_{j}=\mu_{1}\right\}  }.
\end{align*}
The determinant of the Gram matrix for $\left\langle \cdot,\cdot\right\rangle
_{S}$ is an easy consequence by Proposition \ref{inpro1}:
\[
\det\left(  \left\langle h_{\lambda},h_{\mu}\right\rangle _{S}\right)
_{\lambda,\mu\in\mathcal{P}_{n}^{\left(  N\right)  }}=\left(  2^{2n}\left(
\frac{N}{2}+N\kappa\right)  _{2n}\right)  ^{-d\left(  n,N\right)  }%
\widetilde{D}_{n}.
\]
When $N\geq n$ the products in $\widetilde{D}_{n}$ are over all partitions
$\mathcal{P}_{n},\mathcal{P}_{n-1}$ and the inner product formula for
$\left\langle h_{\lambda},h_{\mu}\right\rangle _{h}$ in Proposition
\ref{inpro1} can be considered with $N$ as an indeterminate because
\begin{align*}
m_{\sigma}\left(  1^{N}\right)   &  =m_{\sigma}\left(  1^{l\left(
\sigma\right)  }\right)  \binom{N}{l\left(  \sigma\right)  }\\
&  =m_{\sigma}\left(  1^{l\left(  \sigma\right)  }\right)  \left(  -1\right)
^{l\left(  \sigma\right)  }\left(  -N\right)  _{l\left(  \sigma\right)
}/l\left(  \sigma\right)  !,
\end{align*}
a polynomial in $N$, for $\sigma\in\mathcal{P}$ and $N\geq l\left(
\sigma\right)  $. So the formula for $\widetilde{D}_{n}$ is an identity in the
variables $\kappa,N$ (for $N\geq n$). It would be interesting if there were an
orthogonal basis for $\mathbb{H}_{2n}^{B}$ so that for each element $f$ the
squared norm $\left\langle f,f\right\rangle _{h}$ is a product of integral
powers of factors linear in $\kappa$; this is almost suggested by the nice
form of $\widetilde{D}_{n}$ but so far such a basis has not been found. The
author's suspicion is that the problem is the lack of a sufficiently large set
of commuting self-adjoint operators on $\mathbb{H}_{2n}^{B}$.

\section{Point evaluation of harmonic polynomials}

For almost any classical orthogonal polynomial there is at least one point
where the polynomial can be evaluated in simple terms, for example the Jacobi
polynomials have such values at $\pm1$. It is natural to ask this question
about the symmetric polynomials $h_{\lambda}$. We will consider two special
points with a very simple result for one and a more complicated but quite
interesting analysis for the other.

The easy evaluation is at $x=\left(  1,0,\ldots,0\right)  $ because
$m_{\lambda}\left(  x^{2}\right)  =0$ if $l\left(  \lambda\right)  >1$; thus
for $\lambda\in\widetilde{\mathcal{P}}_{n}^{\left(  N\right)  }$ and $n\geq2 $
we have $h_{\lambda}\left(  x\right)  =\dfrac{2^{-n}B\left(  \lambda,\left(
n\right)  \right)  }{n!\left(  \kappa+\frac{1}{2}\right)  _{n}N}$. By Theorem
\ref{bformula} $B\left(  \lambda,\left(  n\right)  \right)  =\left(
-1\right)  ^{n-\lambda_{1}}n\left(  n-\lambda_{2}-1\right)  !\left(
\prod_{i=2}^{N}\lambda_{i}!\right)  ^{-1}$.

The more interesting evaluation is at $x=1^{N}$ (one reason to single out this
point is that the discriminant $\prod_{i<j}\left(  x_{i}-x_{j}\right)  $
vanishes there of maximum order, another is its role in the generalized
binomial coefficients). Also for $\lambda\in\widetilde{\mathcal{P}}%
_{n}^{\left(  N\right)  }$ the value
\[
h_{\lambda}\left(  1^{N}\right)  =2^{-n}\sum_{\lambda\preceq\mu}\frac{B\left(
\lambda,\mu\right)  }{\mu!\,\left(  \kappa+\frac{1}{2}\right)  _{\mu}}%
\]
depends only on $\lambda$ (and $\kappa$) but not on $N$, so one may as well
assume $N=l\left(  \lambda\right)  $. Analyzing this sum as a rational
function of $\kappa$ we observe that each term is $O\left(  \kappa
^{-n}\right)  $ as $\kappa\rightarrow\infty$. We will prove that $h_{\lambda
}\left(  1^{N}\right)  =O\left(  \kappa^{-\left\lfloor \left(  3n+1\right)
/2\right\rfloor }\right)  $ as $\kappa\rightarrow\infty$; by expanding
$\frac{1}{\left(  \kappa+\frac{1}{2}\right)  _{\mu}}$ as a Laurent series at
$\infty$ and showing that the coefficients of $\kappa^{-i}$ cancel out for
$n\leq i\leq n+\left\lfloor \frac{n-1}{2}\right\rfloor $ (rephrasing this
bound: if $n=2s$ then $2s\leq i\leq3s-1$, if $n=2s+1$ then $2s+1\leq
i\leq3s+1$ where $s=1,2,3,\ldots$). The cancellation is a consequence of the
following property.

\begin{lemma}
\label{bsym0}Suppose $n\geq2,$ $\lambda\in\widetilde{\mathcal{P}}_{n}^{\left(
N\right)  }$ and $g$ is a symmetric polynomial in $N$ variables with
$\deg\left(  g\right)  <n$ then
\[
\sum_{\lambda\preceq\mu}\frac{B\left(  \lambda,\mu\right)  }{\mu!}g\left(
\mu\right)  =0.
\]
\end{lemma}

\begin{proof}
The idea is just a slight generalization of Proposition \ref{bcontig}. Fix
$j<n$ and let $\nu\in\mathcal{P}_{j}^{\left(  N\right)  }$ be arbitrary. We
consider the expansion of $e_{1}^{n-j}m_{\nu}$ in two different ways:
\begin{align*}
e_{1}^{n-j}m_{\nu}  &  =\sum_{\mu\in\mathcal{P}_{n}^{\left(  N\right)  }%
,\nu\subset\mu}%
\genfrac{\langle}{\rangle}{0pt}{}{\mu}{\nu}%
m_{\mu}=\sum_{\nu\subset\mu,\sigma\preceq\mu}B\left(  \sigma,\mu\right)
\genfrac{\langle}{\rangle}{0pt}{}{\mu}{\nu}%
\widetilde{m}_{\sigma}\\
&  =e_{1}^{n-j}\sum_{\tau\preceq\nu}B\left(  \tau,\nu\right)  \widetilde
{m}_{\tau}=\sum_{\tau\preceq\nu}B\left(  \tau,\nu\right)  \widetilde{m}%
_{\tau+\left(  n-j\right)  \varepsilon_{1}}.
\end{align*}
Since $\left\{  \widetilde{m}_{\sigma}:\sigma\in\mathcal{P}_{n}^{\left(
N\right)  }\right\}  $ is a basis it follows that $\sum_{\nu\subset\mu
}B\left(  \lambda,\mu\right)
\genfrac{\langle}{\rangle}{0pt}{}{\mu}{\nu}%
=0$ for any $\lambda\in\widetilde{\mathcal{P}}_{n}^{\left(  N\right)  }$ and
any $\nu\in\mathcal{P}_{j}^{\left(  N\right)  }$ with $0\leq j<n$. The set
$\left\{  \mu\mapsto\mu!%
\genfrac{\langle}{\rangle}{0pt}{}{\mu}{\nu}%
:\nu\in\mathcal{P}_{j}^{\left(  N\right)  },0\leq j<n\right\}  $ is a basis
for the symmetric functions of degree less than $n$. By definition
\begin{align*}
\mu!%
\genfrac{\langle}{\rangle}{0pt}{}{\mu}{\nu}%
&  =\mu!\sum\left\{  \binom{n-j}{\mu-\sigma}:\sigma^{+}=\nu\right\} \\
&  =\sum\left\{  \left(  n-j\right)  !\prod_{i=1}^{N}\frac{\mu_{i}!}{\left(
\mu_{i}-\sigma_{i}\right)  !}:\sigma^{+}=\nu\right\} \\
&  =\sum\left\{  \left(  -1\right)  ^{\left|  \sigma\right|  }\left(
n-j\right)  !\prod_{i=1}^{N}\left(  -\mu_{i}\right)  _{\sigma_{i}}:\sigma
^{+}=\nu\right\}  .
\end{align*}
Note that any terms in the sum with at least one $\sigma_{i}>\mu_{i}$ vanish,
so the sum in the last line can be considered as a symmetric polynomial in the
variables $\mu_{1},\ldots,\mu_{N}$. The highest degree part is clearly
$\left(  n-j\right)  !\,m_{\nu}\left(  \mu\right)  .$ Induction shows that
$m_{\tau}\left(  \mu\right)  \in\mathrm{span}\left\{  \mu!%
\genfrac{\langle}{\rangle}{0pt}{}{\mu}{\nu}%
:\nu\in\mathcal{P}_{j}^{\left(  N\right)  },0\leq j<n\right\}  $ for each
$\tau\in\mathcal{P}_{i}^{\left(  N\right)  },0\leq i<n$.
\end{proof}

The following is a rational version, with an elementary algebraic proof, of
the asymptotic formula for $\Gamma\left(  t\right)  /\Gamma\left(  t+a\right)
$ for $t\rightarrow\infty$.

\begin{lemma}
\label{kfac}There is a sequence $\left\{  p_{j}\left(  n\right)
:j\in\mathbb{N}_{0}\right\}  $ of polynomials such that $p_{j}$ is of degree
$2j$ with leading coefficient $\left(  2^{j}\,j!\right)  ^{-1}$ for $j\geq0$
and
\[
\frac{1}{\left(  \kappa+\frac{1}{2}\right)  _{n}}=\kappa^{-n}\sum
_{j=0}^{\infty}p_{j}\left(  n\right)  \left(  -\kappa\right)  ^{-j}%
\]
for $\left|  \kappa\right|  >n-\frac{1}{2}$, $n\geq0.$
\end{lemma}

\begin{proof}
The coefficients of the Laurent series certainly exist so it suffices to prove
their polynomial property. Clearly $p_{0}\left(  n\right)  =1$ and
$p_{j}\left(  0\right)  =0$ for $j\geq1$. There is a recurrence relation
implicit in the equation:
\begin{align*}
\frac{\kappa}{\left(  \kappa+\frac{1}{2}\right)  _{n+1}}-\frac{1}{\left(
\kappa+\frac{1}{2}\right)  _{n}}  &  =-\frac{n+\frac{1}{2}}{\left(
\kappa+\frac{1}{2}\right)  _{n+1}},\\
\kappa^{-n}\sum_{j=0}^{\infty}\left(  p_{j}\left(  n+1\right)  -p_{j}\left(
n\right)  \right)  \left(  -\kappa\right)  ^{-j}  &  =-\kappa^{-n-1}\sum
_{j=0}^{\infty}p_{j}\left(  n+1\right)  \left(  n+\frac{1}{2}\right)  \left(
-\kappa\right)  ^{-j}.
\end{align*}
Extracting the coefficient of $\left(  -1\right)  ^{j}\kappa^{-n-j}$ in the
second equation yields
\[
p_{j}\left(  n+1\right)  -p_{j}\left(  n\right)  =\left(  n+\frac{1}%
{2}\right)  p_{j-1}\left(  n+1\right)  .
\]
Arguing by induction suppose that $\deg\left(  p_{i}\right)  =2i$ for $0\leq
i\leq j-1$ then the right side of the equation is of degree $2j-1$, the first
difference of $p_{j}$ is of degree $2j-1,$ $p_{j}\left(  0\right)  =0$ and
this determines $p_{j}$ uniquely as a polynomial of degree $2j$. Denote the
leading coefficient of $p_{j}$ by $c_{j}$ then the coefficient of $n^{2j-1}$
in the first difference shows that $2jc_{j}=c_{j-1}$.
\end{proof}

By direct computation $p_{1}\left(  n\right)  =\frac{1}{2}n^{2},\allowbreak
p_{2}\left(  n\right)  =\frac{1}{24}n\left(  n+1\right)  \left(
3n^{2}+n-1\right)  ,\allowbreak p_{3}\left(  n\right)  =\frac{1}{48}%
n^{2}\left(  n+1\right)  \left(  n+2\right)  \left(  n^{2}+n-1\right)  $. The
explanation for the linear factors is that $\frac{1}{\left(  \kappa+_{2}%
^{1}\right)  _{n}}=\left(  \kappa+\frac{1}{2}+n\right)  _{-n}$ for $n<0$,
which is a polynomial of degree $-n$ in $\kappa$, hence $p_{j}\left(
n\right)  =0$ for $j>-n$; the recurrence is also valid for $n<0$.

\begin{proposition}
\label{hat1}For $\lambda\in\widetilde{\mathcal{P}}_{n}^{\left(  N\right)  }$
there is a polynomial $q_{\lambda}\left(  \kappa\right)  $ such that
$\deg\left(  q_{\lambda}\right)  \leq\left\lfloor \frac{n}{2}\right\rfloor
-\lambda_{1}$ and
\[
h_{\lambda}\left(  1^{N}\right)  =\frac{q_{\lambda}\left(  \kappa\right)
}{\left(  \kappa+\frac{1}{2}\right)  _{n}\prod_{i=2}^{N}\left(  \kappa
+\frac{1}{2}\right)  _{\lambda_{i}}}.
\]
\end{proposition}

\begin{proof}
The definition of $\preceq$ implies that
\[
\operatorname{lcm}\left\{  \left(  \kappa+\frac{1}{2}\right)  _{\mu}%
:\lambda\preceq\mu\right\}  =\left(  \kappa+\frac{1}{2}\right)  _{n}%
\prod_{i=2}^{N}\left(  \kappa+\frac{1}{2}\right)  _{\lambda_{i}};
\]
thus the claimed equation is valid for some polynomial $q_{\lambda}$ such that
$\deg\left(  q_{\lambda}\right)  \leq n-\lambda_{1}$ (since $h_{\lambda
}\left(  1^{N}\right)  =O\left(  \kappa^{-n}\right)  $ as $\kappa
\rightarrow\infty$). Assume $\left|  \kappa\right|  >n$ (only a temporary
restriction since we are dealing with an identity of rational functions) and
expand each term in the sum for $h_{\lambda}\left(  1^{N}\right)  $ by Lemma
$\ref{kfac}$ to obtain
\[
h_{\lambda}\left(  1^{N}\right)  =2^{-n}\sum_{\lambda\preceq\mu}\frac{B\left(
\lambda,\mu\right)  }{\mu!}\kappa^{-n}\sum_{\alpha\in\mathbb{N}_{0}^{N}%
}\left(  -\kappa\right)  ^{-\left|  \alpha\right|  }\prod_{i=1}^{N}%
p_{\alpha_{i}}\left(  \mu_{i}\right)  .
\]
For any $s\geq0$ the sum $\sum\limits_{\alpha\in\mathbb{N}_{0}^{N},\left|
\alpha\right|  =s}\left(  -\kappa\right)  ^{-\left|  \alpha\right|  }%
\prod\limits_{i=1}^{N}p_{\alpha_{i}}\left(  \mu_{i}\right)  $ is a symmetric
polynomial of degree $2s$ in $\mu$. In the Laurent series $h_{\lambda}\left(
1^{N}\right)  =\sum_{j=n}^{\infty}a_{j}\kappa^{-j}$ the coefficients $a_{j}=0$
for $2\left(  j-n\right)  <n$ by Lemma \ref{bsym0}; that is $a_{j}=0$ for
$j<\frac{3}{2}n$ and $h_{\lambda}\left(  1^{N}\right)  =O\left(
\kappa^{-\left\lfloor \left(  3n+1\right)  /2\right\rfloor }\right)  $. Thus
$\deg\left(  q_{\lambda}\right)  -\left(  2n-\lambda_{1}\right)  \leq-\left(
n+\left\lfloor \frac{n+1}{2}\right\rfloor \right)  $ and $\deg\left(
q_{\lambda}\right)  \leq\left\lfloor \frac{n}{2}\right\rfloor -\lambda_{1}$.
\end{proof}

The next step is to find an explicit form for $q_{\lambda}$ by imposing
$\left\lfloor \frac{n}{2}\right\rfloor -\lambda_{1}+1$ conditions (the number
of coefficients); these are obtained by a calculation exploiting the fact that
$h_{\lambda}\left(  1^{N}\right)  $ has simple poles at $\kappa=\frac{1}{2}-j$
for $\lambda_{1}+1\leq j\leq n$. There are more points than needed so we
consider only the residues at $\kappa=\frac{1}{2}-j$ for $\left\lfloor
\frac{n+1}{2}\right\rfloor +\lambda_{1}\leq j\leq n$ (note that $\left\lfloor
\frac{n}{2}\right\rfloor +\left\lfloor \frac{n+1}{2}\right\rfloor =n$ so there
are $\left\lfloor \frac{n}{2}\right\rfloor -\lambda_{1}+1$ points in this
list). Only the part $\sum\left\{  \dfrac{B\left(  \lambda,\mu\right)  }%
{\mu!\,\left(  \kappa+\frac{1}{2}\right)  _{\mu}}:\lambda\preceq\mu,\mu
_{1}\geq\left\lfloor \frac{n+1}{2}\right\rfloor +\lambda_{1}\right\}  $ of the
sum for $h_{\lambda}\left(  1^{N}\right)  $ has poles at these points.

Suppose $1\leq j\leq s$ then the residue of $\dfrac{1}{\left(  \kappa+\frac
{1}{2}\right)  _{s}}$ at $\kappa=\frac{1}{2}-j$ equals $\dfrac{\left(
-1\right)  ^{j-1}}{\left(  j-1\right)  !\left(  s-j\right)  !}$ (remove the
factor $\kappa-\frac{1}{2}+j$ and substitute $\kappa=\frac{1}{2}-j$).

\begin{theorem}
For $\lambda\in\widetilde{\mathcal{P}}_{n}^{\left(  N\right)  }$ let $k\left(
\lambda\right)  =\lambda_{1}+\left\lfloor \frac{n+1}{2}\right\rfloor $ then
\begin{align*}
h_{\lambda}\left(  1^{N}\right)   &  =\frac{2^{-n}}{\left(  \kappa+\frac{1}%
{2}\right)  _{k\left(  \lambda\right)  -1}\prod_{i=2}^{N}\left(  \kappa
+\frac{1}{2}\right)  _{\lambda_{i}}}\\
&  \times\sum_{\lambda\preceq\mu,\mu_{1}\geq k\left(  \lambda\right)  }%
\frac{B\left(  \lambda,\mu\right)  }{\mu!}\sum_{j=k\left(  \lambda\right)
}^{\mu_{1}}\frac{\prod_{i=2}^{N}\left(  1-j+\mu_{i}\right)  _{\lambda_{i}%
-\mu_{i}}}{\left(  j-k\left(  \lambda\right)  \right)  !\left(  \mu
_{1}-j\right)  !}\frac{\left(  -1\right)  ^{j-k\left(  \lambda\right)  }%
}{\kappa-\frac{1}{2}+j}.
\end{align*}
\end{theorem}

\begin{proof}
By Proposition \ref{hat1}
\[
h_{\lambda}\left(  1^{N}\right)  =\frac{1}{\left(  \kappa+\frac{1}{2}\right)
_{\lambda+s\varepsilon_{1}}}\sum_{j=k\left(  \lambda\right)  }^{n}\frac{r_{j}%
}{\kappa-\frac{1}{2}+j},
\]
with $s=k\left(  \lambda\right)  -\lambda_{1}-1$ for certain coefficients
$\left\{  r_{j}\right\}  $. Multiply both sides of the equation by $\left(
\kappa+\frac{1}{2}\right)  _{\lambda+s\varepsilon_{1}}$ to obtain
\[
\sum_{j=k\left(  \lambda\right)  }^{n}\frac{r_{j}}{\kappa-\frac{1}{2}%
+j}=2^{-n}\sum_{\lambda\preceq\mu}\frac{B\left(  \lambda,\mu\right)  \left(
\kappa+\frac{1}{2}\right)  _{k\left(  \lambda\right)  -1}}{\mu!\,\left(
\kappa+\frac{1}{2}\right)  _{\mu_{1}}}\prod_{i=2}^{N}\frac{\left(
\kappa+\frac{1}{2}\right)  _{\lambda_{i}}}{\left(  \kappa+\frac{1}{2}\right)
_{\mu_{i}}}.
\]
Note that $\frac{\left(  \kappa+\frac{1}{2}\right)  _{\lambda_{i}}}{\left(
\kappa+\frac{1}{2}\right)  _{\mu_{i}}}=\left(  \kappa+\frac{1}{2}+\mu
_{i}\right)  _{\lambda_{i}-\mu_{i}}$ (for $i\geq2$ and $\lambda\preceq\mu$).
For $j\geq k\left(  \lambda\right)  $ the term corresponding to some $\mu$ has
a pole at $\kappa=\frac{1}{2}-j$ exactly when $j\leq\mu_{1}$ (and
$\lambda\preceq\mu$); the residue of this term is $\dfrac{B\left(  \lambda
,\mu\right)  }{\mu!}\dfrac{\left(  -1\right)  ^{j-1}}{\left(  j-1\right)
!\left(  \mu_{1}-j\right)  !}\left(  1-j\right)  _{k\left(  \lambda\right)
-1}\times\allowbreak\prod\limits_{i=2}^{N}\left(  1-j+\mu_{i}\right)
_{\lambda_{i}-\mu_{i}}$.
\end{proof}

Technically this is not a summation but a transformation formula, but it does
use fewer terms than the original form of $h_{\lambda}\left(  1^{N}\right)  $.
If $\lambda_{1}$ is close to $\frac{n}{2}$ then the formula uses relatively
few terms. There is only one term for $\lambda=\left(  l,l\right)  $ or
$\left(  l,l,1\right)  $ for $n=2l$ or $2l+1$ respectively. It is easy to show
$B\left(  \left(  l,l\right)  ,\left(  2l\right)  \right)  =2\left(
-1\right)  ^{l}$ and $B\left(  \left(  l,l,1\right)  ,\left(  2l+1\right)
\right)  =\left(  -1\right)  ^{l-1}\left(  2l+1\right)  .$ Then the theorem
implies $h_{\left(  l,l\right)  }\left(  1^{N}\right)  =\left(  2^{2l}%
\,l!\,\left(  \kappa+\frac{1}{2}\right)  _{\left(  2l,l\right)  }\right)
^{-1}$ and $h_{\left(  l,l,1\right)  }\left(  1^{N}\right)  =\left(
2^{2l}\,\left(  l-1\right)  !\left(  \kappa+\frac{1}{2}\right)  _{\left(
2l+1,l,1\right)  }\right)  ^{-1}$. The value $h_{\left(  l,l\right)  }\left(
1^{N}\right)  $ can also be found directly by using Kummer's terminating sum
for $_{2}F_{1}\left(  -2l,b;1-2l-b;-1\right)  $ with $b=\frac{1}{2}-\kappa-2l$.

\section{Concluding remarks}

There is an isometry between the space of $\mathbb{Z}_{2}^{N}$-invariant
polynomials with the inner product $\left\langle \cdot,\cdot\right\rangle
_{S}$ and polynomials with the $L^{2}$-inner product for the measure $\left(
\Gamma\left(  \frac{N}{2}+N\kappa\right)  /\Gamma\left(  \kappa+\frac{1}%
{2}\right)  ^{N}\right)  \prod_{i=1}^{N}y_{i}^{\kappa-1/2}dy_{1}\ldots
dy_{N-1},$ on the simplex \newline $\left\{  y\in\mathbb{R}^{N}:\sum_{i=1}%
^{N}y_{i}=1,y_{i}\geq0\text{ each }i\right\}  $, induced by the correspondence
$y=$\newline $\left(  x_{1}^{2},\ldots,x_{N}^{2}\right)  $. Then
$\mathbb{H}_{2n}^{0}$ is isomorphic to the space of polynomials of degree $n$
orthogonal to all polynomials of lower degree$.$ The basis $\left\{
h_{\lambda}:\lambda\in\widetilde{\mathcal{P}}_{n}^{\left(  N\right)
}\right\}  $ maps to a basis for the polynomials symmetric in $\left(
y_{1},\ldots,y_{N}\right)  .$

If one gives up the $B_{N}$-invariance then there is no problem in
constructing a nice orthogonal basis for $\mathbb{H}_{n}$. This basis consists
of products of Jacobi polynomials; a conceptual derivation in terms of
simultaneous eigenfunctions of a set of commuting self-adjoint operators can
be found in \cite[Theorem 2.8]{D3}.

In this paper we constructed a basis for $B_{N}$-invariant spherical harmonics
by introducing a new basis for the symmetric functions. The determinant of the
Gram matrix of the basis was explicitly evaluated. Long ago it was found that
the wave functions of electrons in a crystal with cubical symmetry have
(energy level) degeneracies for $n\geq12$; in our notation $\widetilde
{\mathcal{P}}_{i}^{\left(  3\right)  }=\left\{  \left(  1,1\right)  \right\}
,\left\{  \left(  1,1,1\right)  \right\}  ,\left\{  \left(  2,2\right)
\right\}  ,\left\{  \left(  2,2,1\right)  \right\}  $ for $i=2,3,4,5$
respectively but $\widetilde{\mathcal{P}}_{6}^{\left(  3\right)  }=\left\{
\left(  2,2,2\right)  ,\left(  3,3\right)  \right\}  $. We have not found a
natural and constructive way of orthogonally decomposing $\mathbb{H}_{2n}^{B}
$ when $\dim\mathbb{H}_{2n}^{B}>1$; a candidate for such a decomposition is
the self-adjoint operator $\sum_{1\leq i<j\leq N}\left(  x_{i}\mathcal{D}%
_{j}-x_{j}\mathcal{D}_{i}\right)  ^{4}$ but its eigenvalues are irrational
over the field $\mathbb{Q}\left(  \kappa\right)  $ of rational functions in
$\kappa$. Yet the methods in this paper should give some insight into the
problem of constructing invariant harmonics for the general $B$-type spin
Calogero-Moser model.


\begin{thebibliography}{99}
\bibitem{vD}van Diejen J F 1997 Confluent hypergeometric orthogonal
polynomials related to the rational quantum Calogero system with harmonic
confinement \textit{Commun. Math. Phys.} \textbf{188} 467-497

\bibitem{D1}Dunkl C F 1989 Differential-difference operators associated to
reflection groups \textit{Trans. Amer. Math. Soc.} \textbf{311} 167-183

\bibitem{D2}Dunkl C F 1991 Integral kernels with reflection group invariance
\textit{Canadian J. Math.} \textbf{43} 1213-1227

\bibitem{D3}Dunkl C F 1999 Computing with differential-difference operators
\textit{J. Symbolic Computation }\textbf{28} 819-826

\bibitem{DX}Dunkl C F and Xu Y 2001 Orthogonal Polynomials of Several
Variables \textit{Encyc. Math. Appl. }\textbf{81} (Cambridge: Cambridge
University Press)

\bibitem{K}Kay K G 2001 Exact wave functions from classical orbits: the
isotropic harmonic oscillator and semiclassical applications \textit{Phys.
Rev. A }\textbf{63} 42110

\bibitem{L}Lassalle M 1990 Une formule du bin\^{o}me g\'{e}n\'{e}ralis\'{e}e
pour les polyn\^{o}mes de Jack \textit{C.R. Acad. Sci. Paris S\'{e}r. I Math.}
\textbf{310} 253-256

\bibitem{M}Macdonald I G 1995 \textit{Symmetric Functions and Hall Polynomials
2nd ed.} (Oxford: Clarendon Press)

\bibitem{kT}Taniguchi K 2000 Differential operators that commute with the
$r^{-2}$-type Hamiltonian \textit{Calogero-Moser-Sutherland Models} ed \ J F
van Diejen and L Vinet (New York: Springer) pp 451-459

\bibitem{X1}Xu Y 2000 Harmonic polynomials associated with reflection groups
\textit{Canad. Math. Bull. }\textbf{43} 496-507

\bibitem{YT}Yamamoto T and Tsuchiya O 1996 Integrable $1/r^{2}$ spin chain
with reflecting end \textit{J. Phys A: Math. Gen.}\textbf{\ 29} 3977-3984
\end{thebibliography}
\end{document}